\documentclass{article}
\usepackage[utf8]{inputenc}
\usepackage{epsfig,amssymb}
\usepackage[margin=0.92in]{geometry}
\usepackage{hyperref,diagbox}
\usepackage{graphicx}
\usepackage{epstopdf}
\usepackage{amsmath}
\usepackage{bm}
\usepackage{setspace}
\usepackage{calligra}

\usepackage{mathrsfs}
\usepackage{subcaption,comment}
\usepackage[dvipsnames]{xcolor}
\usepackage{mathtools}
\usepackage{multirow,enumitem}
\usepackage{lipsum}
\usepackage{enumerate}
\usepackage[normalem]{ulem}
\usepackage[affil-it]{authblk}
\usepackage{cite}
\usepackage{cancel}
\usepackage{leftidx,soul}
\usepackage{relsize}
\usepackage{braket}
\usepackage{relsize}
\usepackage{float}

\usepackage[autostyle]{csquotes}

\usepackage[right]{lineno}

\newcommand{\etal}{\textit{et al}.~}
\providecommand{\keywords}[1]{\textbf{\textit{Keywords---}} #1}

\title{Displacement-Driven Approach to Nonlocal Elasticity}
\author{Sansit Patnaik}
\author{Sai Sidhardh}
\author{Fabio Semperlotti}
\affil{Ray W. Herrick Laboratories, School of Mechanical Engineering, Purdue University, West Lafayette, IN 47907}

\begin{document}
\date{}
\maketitle

\begin{abstract}
This study presents a physically consistent displacement-driven reformulation of the concept of action-at-a-distance, which is at the foundation of nonlocal elasticity. In contrast to existing approaches that adopts an integral stress-strain constitutive relation, the displacement-driven approach is predicated on an integral strain-displacement relation. The most remarkable consequence of this reformulation is that the (total) strain energy is guaranteed to be convex and positive-definite without imposing any constraint on the symmetry of the kernels. This feature is critical to enable the application of nonlocal formulations to general continua exhibiting asymmetric interactions; ultimately a manifestation of material heterogeneity. Remarkably, the proposed approach also enables a strong satisfaction of the locality recovery condition and of the laws of thermodynamics, which are not foregone conclusions in most classical nonlocal elasticity theories. Additionally, the formulation is frame-invariant and the nonlocal operator remains physically consistent at boundaries. The study is complemented by a detailed analysis of the dynamic response of the nonlocal continuum and of its intrinsic dispersion leading to the consideration that the choice of nonlocal kernels should depend on the specific material. Examples of exponential or power-law kernels are presented in order to demonstrate the applicability of the method to different classes of nonlocal media. The ability to admit generalized kernels reinforces the generalized nature of the displacement-driven approach over existing integral methodologies, which typically lead to simplified differential models based on exponential kernels. The theoretical formulation is also leveraged to simulate the static response of nonlocal beams and plates illustrating the intrinsic consistency of the approach, which is free from unwanted boundary effects.\\

\noindent\keywords{Nonlocal elasticity, Integral methods, Displacement-driven approach, Energy methods, Nonlocal beams, Nonlocal plates}\\
\noindent All correspondence should be addressed to: \textit{spatnai@purdue.edu} or \textit{fsemperl@purdue.edu}\\

\noindent\textbf{\textit{Highlights}}
\begin{itemize}
    \item Frame-invariant displacement-driven approach to nonlocal elasticity.
    \item Locality recovery condition, boundary effects and thermodynamic consistency analyzed.
    \item Well-posed nonlocal governing equations derived from a convex positive-definite system.
    \item Dispersion relations and material-informed kernels.
    \item Application to the static response of nonlocal beams and plates.
\end{itemize}
\end{abstract}

\section{Introduction}
Experimental and theoretical investigations have shown that, irrespective of the spatial scale, size-dependent or nonlocal effects can become prominent in the response of several structures in the most diverse fields of engineering, nanotechnology, biotechnology, and even medicine \cite{eringen1972linear, bazant1976instability, romanoff2014experimental, romanoff2020review, alotta2020fractional, zhu2020nonlocal, nair2019nonlocal, pradhan2009small}. 
The study and analysis of nonlocal effects has increasingly become an important component of structural design and analysis. Specific examples of size-dependent structures include sensors, biological implants, micro/nano-electromechanical devices and even macroscale structures such as stiffened panels in naval and aerospace applications \cite{zhu2020nonlocal, alotta2020fractional, nair2019nonlocal, pradhan2009small, romanoff2020review}. In the case of micro- and nano-structures, for example those involving carbon nanotubes or graphene sheets, nonlocal effects can be traced back to interatomic interactions, and even medium heterogeneity \cite{pradhan2009small, wang2011mechanisms, trovalusci2017multiscale, tuna2019deformation, patnaik2020towards, lal2019thermomechanical}. These nonlocal forces operate at the nanoscale and hence do not significantly contribute to the overall macroscopic response. For these reasons, nonlocal effects have long been considered negligible at the macro scales. However, at the macro scales, nonlocal effects can originate from material heterogeneities \cite{bazant1976instability, romanoff2014experimental, patnaik2020generalized, trovalusci2014particulate, alotta2020fractional, romanoff2020review, silling2000reformulation, hollkamp2020application} or even from intentionally nonlocal designs \cite{nair2019nonlocal, zhu2020nonlocal}. While long-range interactions occur naturally in heterogeneous materials (e.g. interactions between the layers of a functionally graded material), they can also be induced artificially via carefully designed structural links in the intentional nonlocal designs. 
The existence of size-dependent effects leads to a softening or stiffening of the structure when compared to the predictions made from classical (local) continuum mechanics relations. Additionally, size-dependent structures exhibit anomalous wavenumber-frequency dispersion characteristics leading to wavenumber/frequency dependent wave speeds, contrary to constant wave speeds achieved via local descriptions \cite{fellah2004ultrasonic,fellah2004verification,wang2011mechanisms,buonocore2019occurrence,patnaik2020generalized,lim2015higher}. Consequently, the ability to accurately model size-dependent effects has profound implications for the diverse applications in the aforementioned examples. The inability of the classical continuum theory in capturing size-dependent effects prevented its use in these types of applications and led to the development of nonlocal continuum theories.

Seminal studies by Kroner \cite{kroner1967elasticity}, Edelen \etal \cite{edelen1971thermodynamics}, and Eringen \etal \cite{eringen1972nonlocal,eringen1972linear} laid the theoretical foundation of nonlocal elasticity, and explored its role in the modeling of nonlocal size-dependent structures. The mathematical description of the nonlocal continuum theory proposed in these seminal studies relied on the introduction of additional contributions, resulting from long-range nonlocal interactions, in terms of a convolution integral of the strain field in the constitutive equations. Over the years, several researchers \cite{silling2000reformulation, polizzotto2001nonlocal, polizzotto2004strain, polizzotto2006nonhomogeneous, lazopoulos2006non, askari2008peridynamics, di2013mechanically, pisano2021integral, lazopoulos2020plane, romano2017nonlocal} have proposed different modifications to the formulation presented in the seminal studies \cite{kroner1967elasticity,edelen1971thermodynamics,eringen1972nonlocal}. From a mathematical standpoint, these approaches belong to a class of the so-called \textit{strong integral methods} that capture nonlocal effects by re-defining the stress-strain constitutive law in the form of a convolution integral of either the strain or the stress field over a certain spatial domain (the so-called horizon of nonlocality). 
Depending on whether the nonlocal contributions are modeled using the strain or the stress fields, the integral methods can be classified as strain-driven \cite{eringen1972nonlocal, polizzotto2001nonlocal, polizzotto2004strain, polizzotto2006nonhomogeneous} or stress-driven approaches \cite{romano2017nonlocal}. Further discussions addressing the origin of nonlocal effects, existing theories of nonlocal elasticity, and their applications can be found in this recent review study \cite{shaat2020review}. 

Although the strain- and stress-driven integral approaches to nonlocal elasticity have been able to address a multitude of aspects typical of the response of size-dependent nonlocal structures, some important challenges still remain open. In both these approaches, the nonlocal continuum formulation is based on a nonlocal definition of the stress-strain constitutive relation that adopts a local definition of the strain-displacement kinematic relation. In other terms, the effect of nonlocality is accounted only via the constitutive stress-strain relations \cite{shaat2020review}. Indeed, we will show that it is this assumption that is at the root of several of the challenges faced by existing approaches. These challenges are summarized in the following:
\begin{itemize}
    \item The strain-driven approach leads to a Fredholm integral of the first-kind, which results in ill-posed boundary-value problems. Consequently, its application to model nonlocal structures of practical engineering interest (e.g. beams) leads to inconsistent (also called 'paradoxical') predictions when the same analysis is performed under different loading conditions \cite{challamel2008small, fernandez2016bending, romano2017nonlocal, romano2017constitutive, pisano2021integral}. The issue of an ill-posed formulation was addressed by adopting either stress-driven approaches or enhanced mixed-phase (local and nonlocal) approaches, which led to nonlocal constitutive boundary conditions \cite{romano2017nonlocal, romano2017constitutive, pisano2021integral}. Mathematically speaking, these nonlocal constitutive boundary conditions are additional constraints required to guarantee well-posedness of the integral equations. Notably, the stress-driven approach has yet to achieve a practical numerical implementation (e.g. finite element method) and currently relies on analytical methods to simulate the nonlocal response. While the analytical approach does provide interesting insights, it prevents the application of the stress-driven theory to more complex scenarios involving, as an example, higher dimensional structures subject to general loading conditions or nonlinear behavior.
    
    \item Further, the most direct result of either the strain- or the stress-driven constitutive formulations consists in the inability to yield a positive-definite and convex strain energy density. This follows from the fact that the strain energy density consists of the product between a function and its convolution with a smoothing kernel, which does not yield a quadratic form. Consequently, the only option to achieve a positive-definite strain energy is to require a symmetric nonlocal kernel \cite{polizzotto2001nonlocal, romano2017nonlocal}. While possible, such a condition restricts the application of the resulting theory to structures exhibiting asymmetric long-range interactions \cite{silling2000reformulation,askari2008peridynamics, trovalusci2014particulate, trovalusci2017multiscale,zhu2020nonlocal, hollkamp2020application}. This aspect was also very recently discussed in \cite{batra159misuse} which highlighted the misuse of Eringen's stress-strain nonlocal constitutive relation to model structures made from functionally graded materials.
    
    \item In addition to restrictions on the symmetry of the convolution kernel, several studies also restrict functional nature of the kernel. More specifically, most strain- and stress-driven integral approaches reduce to simplified differential models using the special exponential kernel \cite{pradhan2009small, askes2011gradient, wang2011mechanisms, romano2017nonlocal, lal2019thermomechanical}. The use of an exponential kernel while being a plausible option, prevents the application of the resulting theory to a wide class of structures, particularly those exhibiting anomalous attenuation and dispersion (as we will show in \S\ref{sec: dispersion}). Examples include highly scattering and multifractal media \cite{fellah2004ultrasonic, fellah2004verification, magin2010fractional, hollkamp2019analysis, buonocore2019occurrence, hollkamp2020application} and structures with intentionally designed nonlocal connections \cite{nair2019nonlocal, zhu2020nonlocal}. 
    
    \item The use of a nonlocal stress-strain constitutive relation violates the locality recovery condition which states that, for a physically consistent nonlocal formulation, a uniform stress field and a length-scale independent strain energy must be recovered from a uniform strain field \cite{polizzotto2006nonhomogeneous}. Notably, this limitation was addressed very recently in \cite{pisano2021integral} via an enhanced mixed-phase model, albeit in a weaker form that satisfies only the local stress recovery condition.
    
    \item As noted in \cite{polizzotto2004strain}, the use of a nonlocal stress-strain constitutive behavior is not robust to boundary and surface effects in finite domains. This is a direct result of an inconsistent (and incomplete) truncation of the nonlocal kernel (often, exponential) at boundaries.
\end{itemize}
Some of the above mentioned challenges are also described in \cite{pisano2021integral}. In this study, we attempt to address these conceptual as well as practical limitations of the existing approaches to nonlocal elasticity via fundamental changes in the nonlocal constitutive modeling. More specifically, we reformulate the concept of action-at-a-distance by using nonlocal kinematic relations with a generalized kernel (not restricted to exponential functions). In this study, we will analyze the specific advantages of the proposed approach and illustrate how this method is well equipped to address the different challenges faced by existing integral approaches.

The concept of nonlocal kinematics can be traced back to a seminal study of Drapaca \etal on modeling nonlocal solids by using fractional calculus \cite{drapaca2012fractional}. Fractional-order operators are a class of differ-integral operators, uniquely equipped to model temporal and spatial memory effects in complex materials \cite{fellah2004verification, lazopoulos2006non, magin2010fractional, carpinteri2014nonlocal, mashayekhi2018fractional, patnaik2020generalized, failla2020advanced}. The authors in \cite{drapaca2012fractional} reformulated the classical (local) deformation gradient tensor by using fractional-order derivatives and used the fractional deformation tensor to build the nonlocal continuum theory. This approach was later refined in \cite{sumelka2014thermoelasticity} and applied to the analysis of slender engineering structures. A key limitation of the above mentioned studies \cite{drapaca2012fractional,sumelka2014thermoelasticity} consists in an \textit{ad-hoc} imposition of force and moment balance principles. The governing equations in these studies were derived by a direct replacement of the local stress tensor (in the classical continuum governing equations) with the nonlocal stress tensor. However, as shown in \cite{polizzotto2001nonlocal}, the application of global balance principles to derive governing equations with nonlocal operators leads to the presence of (additional) integral terms. Further, the traction boundary conditions in these studies were identical to Cauchy's postulate for surface tractions. The presence of nonlocal interactions however, leads to a modification of Cauchy's postulate resulting in nonlocal traction conditions \cite{dell2012contact,dell2015postulations,sidhardh2020thermodynamics}. A significant benefit of adopting fractional-order (nonlocal) kinematics was exposed in \cite{patnaik2020ritz, patnaik2020plates, patnaik2020towards}, where the authors could define a positive-definite nonlocal continuum formulation ultimately enabling the use of variational principles. This latter approach had both theoretical advantages including the formulation of well-posed governing equations \cite{patnaik2020ritz, patnaik2020plates, patnaik2020towards}, the incorporation of nonlinear effects \cite{sidhardh2020geometrically, patnaik2020geometrically, sidhardh2020analysis}, and thermodynamic consistency \cite{sidhardh2020thermodynamics}, as well as practical implications like the ability to develop finite element solutions \cite{patnaik2020ritz}. These advantages allowed the application of the fractional-order kinematic approach to model the response of different nonlocal structures subject to mechanical and thermomechanical loading conditions. In all the above mentioned studies, consistent predictions free from boundary and loading effects were obtained.

In this study, we extend the fractional-order kinematic approach proposed in \cite{drapaca2012fractional, sumelka2014thermoelasticity, patnaik2020generalized, patnaik2020towards}, to develop a well-posed and generalized nonlocal kinematic approach. In the following, we will refer to this approach as the \textit{displacement-driven approach to nonlocal elasticity}. Broadly speaking, this approach uses a general differ-integral (nonlocal) definition to capture the strain-displacement relations of a nonlocal continuum. It is exactly this feature that allows the displacement-driven approach to overcome the limitations of existing integral approaches. The process of formulating the displacement-driven approach, analyzing its characteristics, and exploring its application, into three main tasks:
\begin{itemize}
    \item First, the displacement-driven approach to nonlocal elasticity is developed based on generalized nonlocal differ-integral operators introduced within the classical strain-displacement relations. This approach will lead to convex and positive-definite formulation without requiring a symmetric kernel. In this process, we will also analyze the satisfaction of locality recovery condition, the frame-invariance of the formulation, and the consistent behavior of the nonlocal operator at boundaries. We will complete the analysis of the displacement-driven continuum approach by casting the formulation within a thermodynamic framework and demonstrating its thermodynamic consistency. In this regard, we will show how this approach allows a strong (point-wise) imposition of thermodynamic balance laws, consistent with the principles of thermodynamics.
    
    \item Following the development of the displacement-driven constitutive relations, we will derive the equilibrium equations describing the response of the nonlocal solid by using variational principles. Further, we will derive dispersion relations describing the transient dynamics of the nonlocal bulk solid. By taking the example of either an exponential or a power-law kernel, we will demonstrate how the specific choice of the kernel should not be universal but rather determined based on experimental dispersion characteristics carrying the signature of the specific material.
    
    \item Finally, we will apply the displacement-driven approach to model the response of nonlocal beams and plates. For this purpose, we will employ the shear-deformable Timoshenko beam and Mindlin plate models. Since the Euler-Bernoulli beam and the Kirchhoff plate formulations can be recovered as special cases from the Timoshenko and Mindlin formulations, the latter choice allows obtaining more general results. Using a finite element procedure, we will numerically simulate the static response of the nonlocal structural elements for various combinations of nonlocal kernels, kernel parameters, and loading conditions. The choice of using a finite element based solver follows from the flexibility afforded by this method in dealing with general loading conditions as well as integral boundary conditions, typical of nonlocal models \cite{polizzotto2001nonlocal, romano2017nonlocal, patnaik2020towards}. The numerical results will show that the displacement-driven approach predict a consistent (softening) behavior irrespective of the nature of the kernel and of the loading conditions.
\end{itemize}
 
The remainder of the paper is structured as follows: first we discuss the formulation and the relevance of nonlocal elasticity via the displacement driven approach in \S\ref{sec: reformulation}. The analysis of the frame-invariance of the formulation and of the behavior of the nonlocal operator at boundaries follows in \S\ref{sec: FI}. The thermodynamic consistency of the formulation is discussed in \S\ref{sec: thermodynamics}. In \S\ref{sec: GDE}, we derive the governing equations for the displacement-driven approach. The dispersion relations for a 1D unbounded nonlocal solid are derived in \S\ref{sec: dispersion} and are further used to discuss the importance of adopting material-informed kernels. Finally, the proposed approach is applied to model the response of nonlocal beams and plates in \S\ref{sec: application}.

\section{Reformulation of nonlocal elasticity: Displacement-driven approach}
\label{sec: reformulation}
As mentioned in the introduction, the motivation to reformulate the classical approaches to nonlocal elasticity follows from the non-convex nature of the strain energy obtained via these approaches. To better highlight this aspect, we first provide a brief review of the classical integral approaches to nonlocal elasticity. 

In the classical approaches (either strain- or stress-driven) to nonlocal elasticity, the stress and strain at a given point within the nonlocal continuum are related by an integral (convolution) operator. This integral operator enriches the response of a point with the information of the response of a collection of points within a characteristic distance known as the \textit{horizon of nonlocality}. Depending on the specific approach, the stress-strain constitutive relation consists of the convolution of a spatially (monotonically) decaying kernel with either the stress or the strain field in the following fashion \cite{eringen1972nonlocal, polizzotto2001nonlocal, romano2017nonlocal}:
\begin{subequations}
\label{eq: classical_constitutive_expression}
\begin{equation}
    \text{Strain-driven approach: } \sigma_{ij}(\bm{x}) = \int_{\overline{\Omega}} \mathbb{C}_{ijkl} \mathcal{K}_{\varepsilon}(\bm{x},\bm{x}^\prime) \varepsilon_{kl} (\bm{x}^\prime) \mathrm{d}\bm{x}^\prime
\end{equation}
\begin{equation}
    \text{Stress-driven approach: } \varepsilon_{ij}(\bm{x}) = \int_{\overline{\Omega}} \mathbb{S}_{ijkl}\mathcal{K}_{\sigma}(\bm{x},\bm{x}^\prime) \sigma_{kl} (\bm{x}^\prime) \mathrm{d}\bm{x}^\prime
\end{equation}
\end{subequations}
where ${\sigma}_{ij}(\bm{x})$ and ${\varepsilon}_{ij}(\bm{x})$ denote the stress and strain tensor at a point $\bm{x}$ within the domain of the solid $\Omega_\infty$, and $\bm{x}^\prime$ is a point within the horizon of nonlocality $\overline{\Omega} \subseteq \Omega_\infty$. $\bm{\mathbb{C}}$ and $\bm{\mathbb{S}}$ are the fourth-order material elasticity and compliance tensors, respectively. $\mathcal{K}_\varepsilon(\bm{x},\bm{x}^\prime)$ and $\mathcal{K}_\sigma(\bm{x},\bm{x}^\prime)$ are the nonlocal kernels that determine the strength of the long range interactions between the points $\bm{x}$ and $\bm{x}^\prime$, in the strain- and stress-driven approaches, respectively.

The constitutive relations in Eq.~(\ref{eq: classical_constitutive_expression}) result in the following expressions for the total strain energy ($\Pi$) and the strain energy density ($\mathbb{U}(\bm{x})$) of the nonlocal solid \cite{polizzotto2001nonlocal, romano2017nonlocal}:
\begin{subequations}
\label{eq: classical_PE}
\begin{equation}
    \text{Strain-driven approach: } \Pi = \frac{1}{2} \int_{\Omega_\infty} \underbrace{\left[ \varepsilon_{kl} (\bm{x}) \int_{\overline{\Omega}} \mathbb{C}_{ijkl} \mathcal{K}_\varepsilon(\bm{x},\bm{x}^\prime) \varepsilon_{kl} (\bm{x}^\prime) \mathrm{d}\bm{x}^\prime \right]}_{\mathbb{U}(\bm{x}): \text{ Strain energy density}} \mathrm{d}\bm{x}
\end{equation}
\begin{equation}
    \text{Stress-driven approach: } \Pi = \frac{1}{2} \int_{\Omega_\infty} \underbrace{\left[\sigma_{ij}(\bm{x}) \int_{\overline{\Omega}} \mathbb{S}_{ijkl}\mathcal{K}_\sigma(\bm{x},\bm{x}^\prime) \sigma_{kl} (\bm{x}^\prime) \mathrm{d}\bm{x}^\prime \right]}_{\mathbb{U}(\bm{x}): \text{ Strain energy density}} \mathrm{d}\bm{x}
\end{equation}
\end{subequations}
It appears from the above expressions that, in either the strain- or stress-driven approach, the \textit{strain energy density} cannot be reduced to a quadratic form, and hence it is not always guaranteed to be positive-definite. While it is possible to achieve a positive value for the \textit{total strain energy}, this imposes restrictions on the nature of the kernel, that is, the kernel must be positive and symmetric in nature \cite{polizzotto2001nonlocal, romano2017nonlocal}. It was argued, in the introduction, that this latter requirement on symmetry can prevent the application of the resulting theory to several cases of practical relevance \cite{batra159misuse}. As an example, consider a collection of particles that is either heterogeneous or subject to thermal gradients. In such case, the interactions on either side of a given particle are not symmetric \cite{silling2000reformulation,askari2008peridynamics, trovalusci2014particulate, trovalusci2017multiscale, hollkamp2020application, zhu2020nonlocal, batra159misuse}, hence violating the assumption of a symmetric kernel. In other terms, from a physical perspective, the assumption of a symmetric kernel can be considered valid, only for isotropic structures with an underlying symmetry in the arrangement of particles \cite{tuna2019deformation}. 
Mathematically speaking, the restriction of a symmetric kernel is a direct result of the fact that the strain energy density functional is not guaranteed to be positive or convex in nature (simply because it consists of the product of the integral of a function and the function itself, which is not convex in nature). As noted in classical (local) continuum mechanics, a well-posed elastostatic formulation is guaranteed if and only if the strain energy density is both convex and positive-definite in nature.

\begin{figure}[h]
	\centering
	\includegraphics[width=0.8\linewidth]{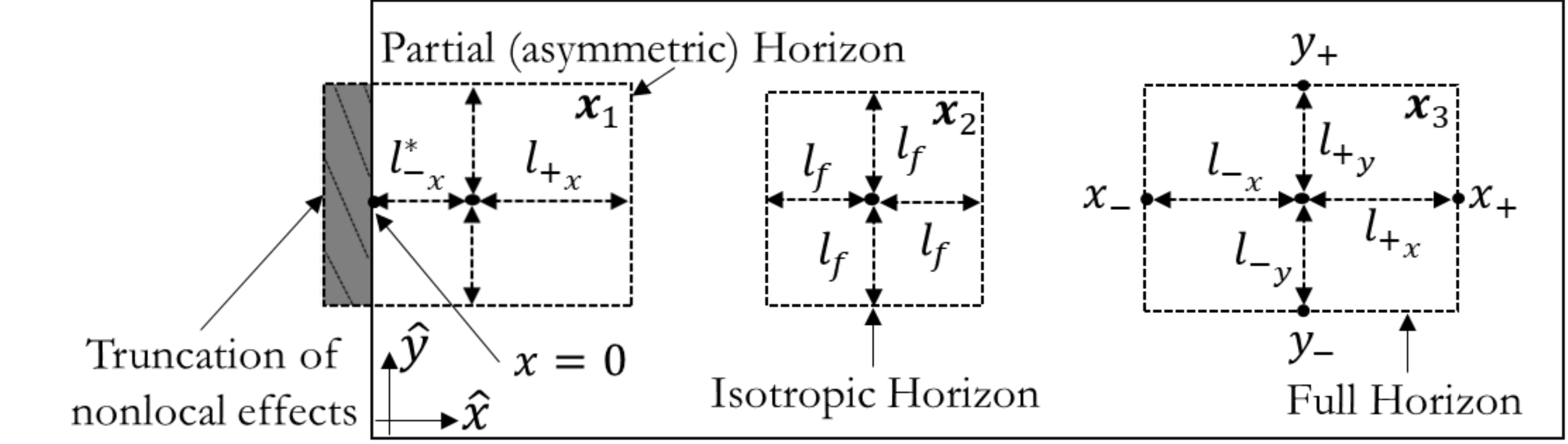}
	\caption{\label{fig: FCM} Schematic illustration of the horizon of nonlocality at different points in a 2D domain. In case of a domain characterized by symmetric nonlocal interactions the horizon of nonlocality at the point $\bm{x}_2$ is equal to $l_f$ on either sides in the $\hat{x}$ direction. However, the horizon of nonlocality at $\bm{x}_1$ must be truncated in the $\hat{x}$, such that $l_{-_x}^*<l_{+_x}$, due to the presence of the boundary. In a similar way, it follows that for a point located on the boundary $x=0$ it must be $l_{-_x}=0$.}
\end{figure}

A possible route to achieve a convex and positive-definite formulation consists in maintaining the quadratic form of the strain energy density, analogous to local elastic continuum mechanics, such that the stress is a work-conjugate of the strain:
\begin{equation}
\label{eq: Displacement_driven_PE}
    \mathbb{U} = \frac{1}{2}\sigma_{ij}\varepsilon_{ij} \equiv \frac{1}{2} \mathbb{C}_{ijkl} \varepsilon_{ij} \varepsilon_{kl} \equiv \frac{1}{2} \mathbb{S}_{ijkl} \sigma_{ij} \sigma_{kl}
\end{equation}
The use of the above quadratic form for the strain energy requires nonlocality to be introduced directly into the strain tensor or equivalently into the stress tensor while using a non-integral (local) stress-strain constitutive relation. Remarkably, this can be achieved by insisting that the kinematic strain-displacement relations are nonlocal (differ-integral) in nature while the stress-strain constitutive relations follow the local continuum mechanics formulation. By following the latter argument, the infinitesimal strain tensor and the Cauchy stress tensor are defined as:
\begin{subequations}
\label{eq: displacemnt_driven_constt_rels}
\begin{equation}
\label{eq: infinitesimal_fractional_strain}
     {\bm{\varepsilon}}=\frac{1}{2} \left[ \overline{\bm\nabla} {\bm{u}} + \overline{\bm\nabla} {\bm{u}}^{T}\right]
\end{equation}
\begin{equation}
\label{eq: stress}
     \bm{\sigma} = \bm{\mathbb{C}}:\bm{\varepsilon}
\end{equation}
\end{subequations}
where $\bm{u}(\bm{x})$ is the displacement field. $\overline{\bm{\nabla}}$ is a differ-integral gradient operator defined as:
\begin{equation}
    \label{eq: grad_operator}
    \overline{\bm{\nabla}} (\cdot) = \overline{D}_{x} (\cdot) \hat{x} + \overline{D}_{y} (\cdot) \hat{y} + \overline{D}_{z} (\cdot) \hat{z}
\end{equation}
where $\{\hat{x},\hat{y},\hat{z}\}$ are the Cartesian basis vectors. $\overline{D}_{x_j}(\cdot)$ is a differ-integral operator similar to the convolution operations in Eq.~(\ref{eq: classical_constitutive_expression}) defined in the following manner:
\begin{equation}
    \label{eq: operator_definition}
    \overline{D}_{x_j}\phi = \int_{\overline{\Omega}_j} c^*_j\mathcal{K}(\bm{x},\bm{x}^\prime) \left[D^1_{x^\prime_j} \phi\right] \mathrm{d}x^\prime_j
\end{equation}
where $\phi$ is an arbitrary function. $c^*_j$ is a scalar multiplier for the differ-integral operator in the $\hat{x}_j$ direction, which will be used to enforce frame-invariance and dimensional consistency of the formulation in \S\ref{sec: FI}. Further, $D^1_{x^\prime_j} \phi$ indicates the classical first-order derivative $\mathrm{d}\phi/\mathrm{d}x^\prime_j$. For a nonlocal horizon $\overline{\Omega}_j = [x_{-_j},x_{+_j}]$ in the $\hat{x}_j$ direction about the point $\bm{x}$, the nonlocal operator can be expanded as: 
\begin{equation}
    \label{eq: operator_definition_2}
    \overline{D}_{x_j}\phi = c_{-_j}^* \int_{x_{-_{j}}}^{x_j} \mathcal{K}(\bm{x},\bm{x}^\prime) \left[D^1_{x^\prime_j} \phi\right] \mathrm{d}x^\prime_j + c_{+_j}^* \int_{x_{j}}^{x_{+_j}} \mathcal{K}(\bm{x},\bm{x}^\prime) \left[D^1_{x^\prime_j} \phi\right] \mathrm{d}{x}^\prime_j 
\end{equation}
where $c_{-_j}^*$ and $c_{+_j}^*$ are independent scalar multipliers. The length of the horizon of nonlocality on either side of the point $\bm{x}$ in the $\hat{x}_j$ direction are denoted as $l_-$ and $l_+$, such that $l_- = x_{j}-x_{-_j}$ and $l_+ = x_{+_j}-x_j$. In general, $l_{-_j} \neq l_{+_j}$ implying that the horizon of nonlocality is not necessarily symmetric like in cases involving material interfaces or boundaries falling within the horizon. This aspect is schematically illustrated in Fig.~(\ref{fig: FCM}).

Revisiting the strain energy of the solid obtained via the above constitutive formulation, we obtain:
\begin{equation}
\label{eq: strain_energy_DD}
    \mathbb{U} = \frac{1}{2} \underbrace{C_{ijkl} \left[ \overline{D}_{x_j} u_i + \overline{D}_{x_i} u_j \right]}_{\text{Nonlocal stress}} \underbrace{\left[ \overline{D}_{x_l} u_k + \overline{D}_{x_k} u_l \right]}_{\text{Nonlocal strain}} 
\end{equation}
which immediately demonstrates its convex nature, irrespective of whether the nonlocal kernel is symmetric or asymmetric in nature. We denominate this nonlocal approach as the \enquote{\textit{displacement-driven approach to nonlocal elasticity}} since nonlocality is embedded via the nonlocal strain-displacement kinematic relations Eq.~(\ref{eq: infinitesimal_fractional_strain}). Note that the stress-strain constitutive relations in Eq.~(\ref{eq: stress}) are non-integral (local) in nature contrary to the classical stress- and strain-driven approaches to nonlocal elasticity \cite{polizzotto2001nonlocal, romano2017nonlocal}. Remarkably, this functional form for the stress-strain constitutive relation trivially satisfies the locality recovery condition \cite{polizzotto2006nonhomogeneous}. More specifically, it ensures that a uniform strain field induces a uniform stress field, and consequently, we obtain a length-scale independent strain energy density via Eq.~(\ref{eq: Displacement_driven_PE}).

The displacement-driven formulation deserves some additional remarks. While the approach was presented starting based on strain energy arguments, the same formulation can be derived via standard continuum mechanics arguments. In this study, we favored the former approach to clearly highlight the connection to a convex and positive-definite nonlocal formulation. When following a continuum mechanics approach, the displacement-driven formulation can be derived by starting from a nonlocal deformation gradient tensor which relates the differential line elements in the deformed and undeformed configurations by using the differ-integral gradient operator defined in Eq.~(\ref{eq: grad_operator}). More specifically, the nonlocal deformation gradient tensor which maps a differential element in the undeformed configuration to the deformed configuration can be defined as:
\begin{equation}
\label{eq: deformation_grad_tensor}
\overline{F}_{ij} = \overline{D}_{X_j} x_i
\end{equation}
where $\bm{X}$ and $\bm{x}$ denote points within the deformed and undeformed coordinates, respectively. Following the above definition for the nonlocal deformation gradient tensor, a fully nonlinear Lagrangian and Eulerian description of the strain tensor can be derived by using the scalar difference of the differential line elements in the deformed and undeformed configurations, similar to classical continuum formulation. The same approach, albeit by using specialized power-law kernels which transforms the gradient operator in Eq.~(\ref{eq: grad_operator}) into a nonlocal gradient operator containing the well established fractional-order derivatives, can be found in \cite{drapaca2012fractional, sumelka2014thermoelasticity, patnaik2020generalized}.

Second, note that the model relies on the important hypothesis of nonlocal kinematics. As highlighted above, in the corresponding continuum description, the differential line elements are defined using differ-integral deformation gradients similar to \cite{drapaca2012fractional, sumelka2014thermoelasticity, patnaik2020generalized}. While this hypothesis might appear unsettling in the context of the more classical continuum approaches, we emphasize that such definition is simply a reformulation of the concept of \textit{action-at-a-distance}. This hypothesis results in assuming that the response of a selected point within the solid is affected directly by the response of a collection of points within the horizon of nonlocality. Given that the differ-integral operator is applied directly to the displacement field (see Eq.~(\ref{eq: infinitesimal_fractional_strain})), from a physical standpoint, the formulation accounts for long-range interactions that are proportional to the relative displacement of distant points within the horizon. It follows that, a change in length of an infinitesimal line at the point $\bm{x}$ between the reference and the current configurations would be affected directly by the response of the points within the nonlocal horizon of $\bm{x}$. This is indeed a possible reformulation of action-at-a-distance that is often implemented in terms of long-range forces \cite{patnaik2020generalized, patnaik2020ritz}.

\section{Frame-invariance: Enforcement and Implications}
\label{sec: FI}
The reformulation of the kinematic relations via differ-integral operators requires a thorough investigation of the frame-invariance, which states that a rigid-body motion should not induce a change in the strain energy of the body. In this section, we will derive the necessary constraints for the scalar multipliers $c_-^*$ and $c_+^*$ such that frame-invariance is satisfied at all points within the nonlocal continuum. Remarkably, we will show that the satisfaction of frame-invariance via the scalar multipliers also ensures a consistent behavior of the nonlocal operator at material boundaries. 

\subsection{Enforcing frame-invariance}
Consider a rigid-body motion superimposed on the reference configuration of the body as:
\begin{equation}
\label{sq5}
\bm{\chi}(\bm{X},t)=\bm{C}(t)+\bm{Q}(t)\bm{X},
\end{equation}
where $\bm{C}(t)$ is a spatially constant term representing a translation and $\bm{Q}(t)$ is a proper orthogonal tensor denoting a rotation. Under this rigid-body motion, the deformation gradient tensor $\overline{\bm{F}}^{\bm{\chi}}$ should be an orthogonal tensor such that ${\overline{\bm{F}}^{\bm{\chi}}}^T \overline{\bm{F}}^{\bm{\chi}} = \textbf{I}$. More specifically, the nonlocal deformation gradient tensor should transform as $\overline{\bm{F}}^{\bm{\chi}} = \bm{Q}$ in order to ensure that the strain measures corresponding to the rigid-body motion are null. From the definition of $\overline{\bm{F}}^{\bm{\chi}}$ given in Eq.~(\ref{eq: deformation_grad_tensor}) it follows that:
\begin{equation}
\label{sq6}
\begin{split}
\overline{D}_{x_j}\chi_i =  c_{-_j}^* \int_{x_{-_{j}}}^{x_j} \mathcal{K}(\bm{x},\bm{x}^\prime) {D^1_{x^\prime_j} \chi_{i}(\bm{x}^\prime,t)} \mathrm{d}x^\prime_j + c_{+_j}^* \int_{x_{j}}^{x_{+_j}} \mathcal{K}(\bm{x},\bm{x}^\prime) {D^1_{x^\prime_j} \chi_{i}(\bm{x}^\prime,t)} \mathrm{d}{x}^\prime_j 
\end{split}
\end{equation}
Further, noting that $D^1_{x^\prime_j} C_i(t)=0$ and $\bm{Q}=\bm{Q}(t)$ we obtain:
\begin{equation}
\label{sq8}
D^1_{x^\prime_j} \chi_{i}(\bm{x}^\prime,t) = Q_{ik}x_{k,j} = Q_{ik}\delta_{kj} = Q_{ij}
\end{equation}
Thus, under the rigid body motion $\bm{\chi}$:
\begin{equation}
\label{sq9}
\overline{D}_{x_j} \chi_i = \left[ c_{-_j}^* \int_{x_{-_{j}}}^{x_j} \mathcal{K}(\bm{x},\bm{x}^\prime) \mathrm{d}x^\prime_j + c_{+_j}^* \int_{x_{j}}^{x_{+_j}} \mathcal{K}(\bm{x},\bm{x}^\prime) \mathrm{d}\bm{x}^\prime_j \right]Q_{ij}
\end{equation}
To achieve frame-invariance, that is enforcing $\overline{D}_{x_j} \chi_i = Q_{ij}$ at all points within the continuum and at all time instants, we require that:
\begin{subequations}
\label{eq: contants_expression}
\begin{equation}
c_{-_j}^* = \frac{1}{2\int_{x_{-_{j}}}^{x_j} \mathcal{K}(\bm{x},\bm{x}^\prime) \mathrm{d}x^\prime_j}
\end{equation}
\begin{equation}
c_{+_j}^* = \frac{1}{2\int_{x_{j}}^{x_{+_j}} \mathcal{K}(\bm{x},\bm{x}^\prime) \mathrm{d}\bm{x}^\prime_j}
\end{equation}
\end{subequations}
Notably, the above expressions for $c_{-_j}^*$ and $c_{+_j}^*$ also ensure the dimensional consistency of the formulation. Here below we provide the expressions of $\{c_{-_j}^*, c_{+_j}^*\}$ for some nonlocal kernels commonly adopted in literature:
\begin{itemize}
    \item \textit{Exponential kernel}: $\mathcal{K}(\bm{x},\bm{x}^\prime) = \exp \left( -\frac{|\bm{x} - \bm{x}^\prime|}{l_0} \right)$, where $l_0$ is a material parameter:
    \begin{subequations}
    \label{eq: constants_exp_kernel}
    \begin{equation}
    c_{-_j}^* = \frac{1}{2 l_0} \left[ 1 - \exp^{\left(-\frac{l_-}{l_0} \right)} \right]^{-1}
    \end{equation}
    \begin{equation}
    c_{+_j}^* = \frac{1}{2 l_0} \left[ 1 - \exp^{\left(-\frac{l_+}{l_0} \right)} \right]^{-1}
    \end{equation}
    \end{subequations}
    
    \item \textit{Power-law kernel}: $\mathcal{K}(\bm{x},\bm{x}^\prime) = |\bm{x} - \bm{x}^\prime|^{-\alpha}/\Gamma(1-\alpha)$, where $\Gamma(\cdot)$ is the well known Gamma function and $\alpha \in (0,1]$ is a material parameter:
    \begin{subequations}
    \label{eq: constants_pl_kernel}
    \begin{equation}
    c_{-_j}^* = \frac{1}{2} \Gamma(2-\alpha) l_-^{\alpha-1}
    \end{equation}
    \begin{equation}
    c_{+_j}^* = \frac{1}{2} \Gamma(2-\alpha) l_+^{\alpha-1}
    \end{equation}
    \end{subequations}
\end{itemize}
While the above results are limited to the exponential and the power-law kernels, the approach is general and can be extended to other nonlocal kernels such as, for example, bell and conical curves employed in classical strain- and stress-driven models of nonlocal elasticity \cite{polizzotto2001nonlocal}. In each case, the constants $c_-^*$ and $c_+^*$ must be derived by using Eq.~(\ref{eq: contants_expression}). We will analyze the specific impact of different nonlocal kernels in more detail in \S\ref{sec: dispersion}.

\subsection{Behavior of the nonlocal operator at boundaries}
A detailed analysis of the expressions for the scalar multipliers $c_{-_j}^*$ and $c_{+_j}^*$ in Eq.~(\ref{eq: contants_expression}) presents some interesting insights on the nature of the kernel at material boundaries that will be of practical relevance. For the most general definition of the kernel $\mathcal{K}(\bm{x},\bm{x}^\prime)$ and the location of the point $\bm{x}$ recall that $l_-\neq l_+$. It immediately follows from Eq.~(\ref{eq: contants_expression}) $c_{-_j}^* \neq c_{+_j}^*$. This latter condition, apart from ensuring an exact satisfaction of frame-invariance, also ensures the completeness of the nonlocal kernel at points close to the material boundaries or at points located on the boundary itself. In order to demonstrate this latter statement, we investigate the behavior of the nonlocal operator $\overline{D}_{x_j} \phi$ at a point located on the boundary of the nonlocal solid given in Fig.~(\ref{fig: FCM}).  
    
For a point (say $\bm{x}_0$) located on one of the boundaries (identified by the normal in the $j^{th}$ direction), either $l_{-_j} \rightarrow 0$ or $l_{+_j} \rightarrow 0$ (see Fig.~(\ref{fig: FCM})). Here below we present the expression for the nonlocal operator when $l_{-_j}\rightarrow 0$ (similar expressions hold when $l_{+_j} \rightarrow 0$). This limiting case, upon using the expression for $c_{-_j}^*$ in Eq.~(\ref{eq: contants_expression}), gives:
\begin{equation}
\label{sq10}
\lim_{l_{-_j} \to 0} \overline{D}_{x_j} \phi = \lim_{l_{-_j} \to 0} \left[ \frac{1}{2\int_{x_{0_j} - l_{-_j}}^{x_{0_j}} \mathcal{K}(\bm{x},\bm{x}^\prime) \mathrm{d}x^\prime_j}  \int_{x_{0_j} - l_{-_j}}^{x_{0_j}} \mathcal{K}(\bm{x},\bm{x}^\prime) D^1_{x^\prime_j} \phi \mathrm{d}x^\prime_j + c^*_{+_j} \int_{x_{0_j}}^{x_{+_j}} \mathcal{K}(\bm{x},\bm{x}^\prime) D^1_{x^\prime_j} \phi \mathrm{d}x^\prime_j \right]
\end{equation}
Since the interval length ($=l_{-_j}$) of the left integral within the right-hand side of the above equation is very small, $D^1_{x^\prime_j} \phi$ can be assumed to be constant and equal to the boundary condition:
\begin{equation}
\label{sq11}
D^1_{x^\prime_j} \phi = D^1_{x_j} \phi \Big|_{\bm{x}_0}
\end{equation}
Substituting Eq.~(\ref{sq11}) in Eq.~(\ref{sq10}) leads to the following:
\begin{equation}
\label{sq12}
\lim_{l_{-_j} \to 0} \overline{D}_{x_j} \phi = \lim_{l_{-_j} \to 0} \left[ \frac{D^1_{x_j} \phi \Big|_{\bm{x}_0}}{2 \cancelto{}{\int_{x_{0_j} - l_{-_j}}^{x_{0_j}} \mathcal{K}(\bm{x},\bm{x}^\prime) \mathrm{d}x^\prime_j}}  \cancelto{}{\int_{x_{0_j} - l_{-_j}}^{x_{0_j}} \mathcal{K}(\bm{x},\bm{x}^\prime) \mathrm{d}x^\prime_j} + c^*_{+_j} \int_{x_{0_j}}^{x_{+_j}} \mathcal{K}(\bm{x},\bm{x}^\prime) D^1_{x^\prime_j} \phi \mathrm{d}x^\prime_j \right]
\end{equation}
which immediately yields:
\begin{equation}
\label{sq13}
\lim_{l_{-_j} \to 0} \overline{D}_{x_j} \phi = \frac{1}{2} \underbrace{D^1_{x_j} \phi \Big|_{\bm{x}_0}}_{\substack{\text{Local effect due} \\ \text{to truncation of} \\ \text{nonlocal horizon}}} + \underbrace{c^*_{+_j} \int_{x_{0_j}}^{x_{+_j}} \mathcal{K}(\bm{x},\bm{x}^\prime) D^1_{x^\prime_j} \phi \mathrm{d}x^\prime_j}_{\text{Remaining nonlocal interaction}}
\end{equation}
From Eq.~(\ref{sq13}), it is immediate to observe that while the right-handed operator captures nonlocality ahead of the point $\bm{x}_0$ (in the $j^{th}$ direction), the left-handed operator is reduced to the classical first-order derivative. From a mathematical standpoint, the definition of the nonlocal operator ensures that, independent of the nature of the kernel, we recover a delta-function (and a corresponding local behavior) every time the size of the nonlocal horizon on any given side of a point approaches (or is equal to) zero. This suggests that the truncation of the nonlocal horizon and the corresponding convolution at the boundary has been accounted for in a consistent manner. From an atomistic or molecular dynamic perspective, the reduction of the left-handed operator to the classical first-order derivative is physically analogous to the annihilation (or removal) of long-range interactions between atoms or molecules close to material boundaries \cite{tuna2019deformation}.

Before proceeding further, we merely note that the equality $c_{-_j}^* = c_{+_j}^*$ holds if and only if the nonlocal kernel is symmetric and $l_{-_j} = l_{+_j}$. The above conditions require that the body has uniform material properties (such that the nonlocal kernel exhibits uniform spatial characteristics) and the point under consideration has a symmetric horizon of nonlocality (possible for points located sufficiently far from the material boundaries).

\section{Considerations on the thermodynamics}
\label{sec: thermodynamics}
In this section, we explore the implications of a displacement-driven approach on the thermodynamic formulation of the nonlocal solid and we will show how this approach is uniquely equipped to enforce the laws of thermodynamics in a strong (point-wise) sense.\\

\noindent\textit{First law of thermodynamics:} For a nonlocal solid, the statement of the conservation of energy for a point $\bm{x}$ includes, in addition to the local strain energy, the energy associated with the long-range interactions with points contained within its horizon of nonlocality. The functional relationship between the energy density and the different energy exchanges (local and nonlocal) is typically expressed as $e=e(\bm{\varepsilon}_l,\mathcal{R}(\bm{\varepsilon}_l),{\eta})$ where $\bm{\varepsilon}_l$ is the local strain field, $\mathcal{R}(\bm{\epsilon})$ denotes a linear integral operator which models nonlocality in the solid, and $\eta$ is the entropy of the solid. Broadly speaking, the nonlocal strain in the displacement driven approach allows combining the local strain $\bm{\varepsilon}_l$ and its integral $\mathcal{R}(\bm{\varepsilon}_l)$. More specifically, given the nonlocal kinematic relations described in \S\ref{sec: reformulation}, the contribution of the energy contained in the long-range interactions are fully captured in the nonlocal strain. In this context, recall the argument that the nonlocal kinematic relations were a reformulation of the concept of action-at-a-distance, and hence carry information of nonlocal interactions within the solid. It immediately follows that the internal energy density of the solid can be expressed explicitly as a function of the nonlocal strain. Mathematically, we have $e=e(\bm{\varepsilon},\eta)$ $\bm{\varepsilon}$ is the nonlocal strain and $\eta$ is the entropy of the solid. This latter conclusion is remarkable as it allows the first law of thermodynamics to be applied in a strict sense at every point in the domain:
\begin{equation}
    \label{eq: first_law_thermod}
    \dot{e} = {\sigma}_{ij}{\varepsilon}_{ij} + h - q_{i,i}~\forall \bm{x}\in\Omega
\end{equation}
where the notation $(\dot{\square})$ denotes the first-order derivative with respect to time, $h$ is the heat generated internally per unit volume, and $\bm{q}$ is the heat flux density.\\

\noindent\textit{Second law of thermodynamics:} The second law of thermodynamics states that the internal entropy production rate within a solid is non-negative for all points inside the solid, that is $\dot{\eta}_{0} \geq 0~\forall~\bm{x}\in\Omega$. In order to apply the second law of thermodynamics to the displacement-driven formulation, we consider the internal entropy production rate:
\begin{equation}
    \label{eq: entropy_int_gen}
    \dot{{\eta}}_{\text{0}}=\dot{{\eta}}-\left[\frac{h}{T}-\left(\frac{q_i}{T}\right)_{,i}\right]
\end{equation}
where $\nabla(\cdot)$ denotes the classical divergence operator and $T$ denotes the absolute temperature of the solid. In analogy with the classical approach, we introduce the Legendre transformation $\psi = e - T \eta$, where $\psi$ denotes the Helmholtz free energy. By following the rationale presented in defining $e=e(\varepsilon,\eta)$, it immediately follows that $\psi = \psi(\varepsilon,T)$. By using the Legendre transformation along with Eq.~(\ref{eq: entropy_int_gen}) we obtain:
\begin{equation}
    \label{eq: clausius_duhem}
    T\dot{{\eta}}_{\text{0}} = {\sigma}_{ij}{\varepsilon}_{ij} - \dot{\psi} - \eta \dot{T} - T_{,i}\frac{q_{i}}{T} \geq 0
\end{equation}
Note that the statement of the second law of thermodynamics for the displacement-driven approach is analogous to the classical Clausius-Duhem inequality in its function form, hence presenting a clear departure from several classical nonlocal approaches characterized by the presence of additional terms. These additional terms, resulting from a functional dependence of $\psi$ on $\mathcal{R}(\bm{\varepsilon}_l)$, casting the second law of thermodynamics in strong-form. Indeed, these terms disappear only when a weak-form is considered (that is $\int_{\Omega}\dot{{\eta}}_{\text{0}}\geq 0$). The ability to satisfy the second law of thermodynamics only in a weak sense is the main source of physical inconsistencies of classical nonlocal formulation, as highlighted in \cite{polizzotto2001nonlocal}.

The Clausius-Duhem inequality in Eq.~(\ref{eq: clausius_duhem}) can be used to derive thermodynamically consistent constitutive equations which should be equal to the ones presented in \S\ref{sec: reformulation}. By substituting the functional relationship $\psi = \psi(\varepsilon,T)$, the inequality in Eq.~(\ref{eq: clausius_duhem}) is simplified as:
\begin{equation}
\label{eq: constt_1}
   T\dot{{\eta}}_{\text{0}} = \left({\sigma}_{ij}-\frac{\partial \psi}{\partial {\varepsilon}_{ij}}\right)\dot{{\varepsilon}}_{ij}-\left(\eta+\frac{\partial \psi}{\partial T}\right)\dot{T}-T_{,i}\frac{q_{i}}{T}\geq 0
\end{equation}
Since the above inequality must hold for arbitrary choices of the independent fields ${\bm{\varepsilon}}$ and $T$ at all times, we obtain the following constitutive laws:
\begin{subequations}
\label{constt_eq}
\begin{equation}
\label{eq: mech_constt_eq}
    {\sigma}_{ij}=\frac{\partial \psi}{\partial {\varepsilon}_{ij}} ~\forall\bm{x} \in \Omega
\end{equation}
\begin{equation}
\label{eq: therm_constt_eq}
    {\eta}=-\frac{\partial \psi}{\partial T} ~\forall\bm{x} \in \Omega
\end{equation}
\end{subequations}
Finally, by using the above constitutive relations within Eq.~(\ref{eq: clausius_duhem}), we obtain:
\begin{equation}
\label{eq: reduced_ineq}
    T\dot{{\eta}}_{\text{int}}=-T_{,i}\frac{q_{i}}{T}\geq 0 ~\forall\bm{x} \in \Omega
\end{equation}
which establishes the strong-form of the second law of thermodynamics.

\section{Governing equations}
\label{sec: GDE}
In this section, we derive the strong-form of the equilibrium equations governing the response of the nonlocal solid by using variational principles. Note that the governing equations can also be derived by using a more classical Newtonian approach based on force and moment balance over a representative volume element of the nonlocal domain. Although both approaches yield identical results, the variational approach is deemed more direct given that the explicit use of forces and moments balance would require additional care in accounting for the effect of the nonlocal interactions, as described in detail in \cite{polizzotto2001nonlocal, sidhardh2020thermodynamics}. 

We derive the strong-form of the governing equations by using the Hamilton's principle:
\begin{equation}
    \label{eq: extended_hamiltons_principle1}
    \int_{t_1}^{t_2} \delta \left(\Pi -  V -  T\right)\mathrm{d}t=0
\end{equation}
In the above equation, $V$ and $T$ denote the work done by externally applied forces and the kinetic energy of the system. Recall that $\Pi$ denotes the total strain energy of the formulation. Substituting the expressions for the different physical variables, we obtain:
\begin{equation}
    \label{eq: extended_hamiltons_principle2}
    \int_{t_1}^{t_2} \delta\Bigg[ \underbrace{ \frac{1}{2} \int_\Omega (\bm\sigma:\bm\varepsilon) \mathrm{d}\mathbb{V}}_{\text{Strain Energy}} - \underbrace{\int_\Omega (\bm{\overline{b}}\cdot\bm{u}) \mathrm{d}\mathbb{V}  - \int_{\partial \Omega} (\bm{\overline{t}}\cdot\bm{u}) \mathrm{d}\mathbb{A}}_{\text{External Work}} -  \underbrace{\frac{1}{2} \int_{\Omega} \rho (\bm{\dot{u}} \cdot {\bm{\dot{u}}})\mathrm{d}\mathbb{V}}_{\text{Kinetic Energy}} \Bigg] \mathrm{d}t=0
\end{equation}
where $\mathrm{d}\mathbb{V}$ and $\mathrm{d}\mathbb{A}$ indicate volume and area elements of the nonlocal solid, respectively. $\rho$ denotes the density of the solid, and $\bm{\overline{b}}$ and $\bm{\overline{t}}$ are prescribed values of the body force per unit volume and the surface traction per unit area, respectively. 

By simplifying Eq.~(\ref{eq: extended_hamiltons_principle2}) using principles of variational calculus, we obtain the governing equation as:
\begin{equation}
\label{eq: GDE}
\bm{\widetilde{\nabla}}\cdot\bm{\sigma} + \bm{\overline{b}} = \rho\ddot{\bm{u}}
\end{equation}
The associated boundary conditions are obtained as:
\begin{equation}
\label{eq: BCs}
\bm{\tilde{I}}_{\bm{\hat{n}}} \cdot \bm{\sigma} = \bm{\overline{t}}  ~~~~ \text{or} ~~~~ {\bm{u}} = \bm{\overline{u}} ~~~~ \forall ~ \bm{x} \in \partial\Omega
\end{equation}
The operator $\bm{\tilde{I}}_{\bm{\hat{n}}}(\cdot)$ is defined as:
\begin{equation}
\label{eq: integral_normal_operator}
    \bm{\tilde{I}}_{\bm{\hat{n}}}(\cdot) = n_x \tilde{I}_{x}(\cdot) \hat{x} + n_y \tilde{I}_{y}(\cdot) \hat{y} + n_z \tilde{I}_{z}(\cdot) \hat{z}
\end{equation}
where $\bm{\hat{n}}=n_x \hat{x} + n_y \hat {y} + n_z \hat{z}$ is the unit normal vector at the surface of the solid. In the above equation, $\tilde{I}_{x_j}(\cdot)$ is an integral operator defined in the following manner: 
\begin{equation}
\label{eq: integral_def}
    \tilde{I}_{x_j} \phi = c_{-_j}^* \int_{x_j-l_{+_j}}^{x_j} \mathcal{K}(\bm{x},\bm{x}^\prime) \phi \mathrm{d}x_j^\prime + c_{+_j}^* \int^{x_j + l_{-_j}}_{x_j} \mathcal{K}(\bm{x},\bm{x}^\prime) \phi \mathrm{d}x_j^\prime
\end{equation}
Further, the gradient operator denoted by $\bm{\widetilde{\nabla}}(\cdot)$ is a differ-integral operator (analogous to $\overline{\nabla}$ in Eq.~(\ref{eq: operator_definition})) containing the above defined integral operators:
\begin{equation}
    \label{eq: modified_grad_operator}
    \bm{\widetilde{\nabla}}^{\alpha_m} (\cdot) = D^1_x \left[ \tilde{I}_{x}(\cdot) \right] \hat{x} + D^1_y \left[ \tilde{I}_{y}(\cdot) \right] \hat{y} + D^1_z \left[ \tilde{I}_{z}(\cdot) \right] \hat{z} \equiv \widetilde{D}_{x} (\cdot) \hat{x} + \widetilde{D}_{y} (\cdot) \hat{y} + \widetilde{D}_{z} (\cdot) \hat{z}
\end{equation}
The detailed derivation of the governing equations is provided in the Appendix. Note that the natural boundary conditions are integral (and hence, nonlocal) in nature, consistent with existing integral approaches to nonlocal elasticity \cite{polizzotto2001nonlocal, romano2017nonlocal}. From a physical perspective, the nonlocal natural boundary conditions account for the combined effect of long-range forces exerted on the boundary points by points located within the nonlocal solid. While the use of variational principles allowed obtaining the nonlocal natural boundary conditions in a relatively straight-forward manner, the use of global conservation principles (mechanical balance laws) requires a careful modification of the surface tractions to account for the nonlocal interactions of the boundary points (mentioned previously). In other terms, Cauchy's postulate for surface tractions requires a modification to account for the nonlocal interactions. This latter aspect has been addressed in detail in \cite{dell2012contact, dell2015postulations, sidhardh2020thermodynamics}. 

\section{Dispersion relations and choice of kernels}
\label{sec: dispersion}
The nonlocal elastodynamic framework presented in \S\ref{sec: GDE} allows for the derivation of wavenumber-frequency dispersion relations describing the transient response of the nonlocal bulk solid. It is anticipated that the dispersion relations functionally depend on the nature of the nonlocal kernel used to model nonlocality within the solid. In this section, we will derive the dispersion relations for an unbounded 1D nonlocal solid and use them to gain insight into the specific effects of the nature of the nonlocal kernel. 

The governing equation describing a 1D nonlocal solid can be extracted from Eq.~(\ref{eq: GDE}) as:
\begin{equation}
\label{eq: 1D_GDE}
    E \widetilde{D}_x \overline{D}_x u = \rho \ddot{u}
\end{equation}
where $E$ is the modulus of elasticity of the isotropic solid. To obtain the dispersion relation, we substitute in the elastodynamic equation the following ansatz:
\begin{equation}
    \label{eq: 1D_solution_assumption}
    u(x,t) = u_0 e^{i(kx - \omega t)}
\end{equation}
where $u_0$ is the amplitude of the longitudinal wave, $k$ denotes the wave-number, $\omega$ denotes the angular frequency of free longitudinal vibrations, and $i=\sqrt{-1}$. Using the expressions for the differ-integral operators in Eq.~(\ref{eq: 1D_GDE}), we obtain the following expression for the dispersion relation in an unbounded nonlocal solid:
\begin{equation}
    \label{eq: dispersion_step_1}
    \omega^2 + ik e^{-ikx} \frac{E}{\rho} D^1_x \left[ \int^\infty_{-\infty} \mathcal{K}(x,x^\prime) \int^\infty_{-\infty} \mathcal{K}(x^\prime, s^\prime) e^{iks^\prime} \mathrm{d} s^\prime  \mathrm{d}x^\prime \right] = 0
\end{equation}
The above expression provides important insights into the nature of the formulation. First, the expression suggests that for a well-posed formulation resulting in bounded wave speeds, it is essential that the convolution integrals in Eq.~(\ref{eq: dispersion_step_1}) exist and are bounded. This latter condition is always satisfied if the nonlocal kernel: 1) decays spatially with the increasing inter-particle distance, that is, $D^1_{x-x^\prime}\mathcal{K}(x,x^\prime) < 0 ~\forall \{x,x^\prime\}$, and 2) is positive-definite (more concretely, the Fourier transform of the kernel is positive everywhere \cite{bavzant1984instability}). Note also that, a decaying kernel, also guarantees the existence of the integrals in Eq.~(\ref{eq: contants_expression}) which are essential to enable a frame-invariant formulation. This condition is also consistent with existing approaches to nonlocal elasticity, which invariably use decaying kernels \cite{eringen1972nonlocal, polizzotto2001nonlocal, lazopoulos2006non, askes2011gradient, romano2017nonlocal}. 
At this stage, following the analysis of frame-invariance and of the dispersion relations, we are well-equipped to summarize all the sufficient requirements for the nonlocal kernel in order to ensure a physically consistent formulation: the kernel should be positive-definite and must decay monotonically with the interparticle distance.

Note that further simplifications of the expression in Eq.~(\ref{eq: dispersion_step_1}) requires the specific functional definition of the nonlocal kernel. From a practical perspective, this suggests that the applicability of a given functional form of the nonlocal kernel, to model a certain class of nonlocal structures, depends on whether the dispersion characteristics obtained by simplification of Eq.~(\ref{eq: dispersion_step_1}) match the experimentally obtained dispersion characteristics of the class of materials under consideration. This latter aspect becomes evident by considering the example of the exponential and the power-law kernels introduced in \S\ref{sec: reformulation}. The dispersion relations for the exponential and the power-law kernels are obtained by using Eq.~(\ref{eq: dispersion_step_1}) as:
\begin{subequations}
\label{eq: dispersion_examples}
\begin{equation}
\label{eq: dispersion_exp}
    \text{Exponential kernel: } \left(\frac{\omega}{k}\right)^2 = {\frac{E}{\rho}} \left(\frac{1}{1+k^2l_0^2} \right)^2
\end{equation}
\begin{equation}
\label{eq: dispersion_pl}
    \text{Power-law kernel: } \left(\frac{\omega}{k}\right)^2 =  {\frac{E}{\rho}} \big[ \underbrace{\cos \left( {\pi+\alpha\pi} \right)}_{\text{Propagation}} + i \underbrace{\sin \left( {\pi+\alpha\pi} \right)}_{\text{Attenuation}} \big] \left( kl^* \right)^{2(\alpha-1)}
\end{equation}
\end{subequations}
While the derivation of the dispersion relation is straightforward for the exponential kernel (given the exponential nature of both the kernel and the assumed wave solution), some special relations are required to derive the dispersion relations for the power-law kernel. These relations (that concern the fractional-order derivatives of an exponential function) as well as the derivation of the dispersion relations for the power-law kernel can be found in \cite{patnaik2020towards}. In the power-law dispersion, $l^*=1$m is a dimensional constant which is used to enforce dimensional consistency for the 1D unbounded solid \cite{hollkamp2019analysis,patnaik2020towards}. This constant is required for the 1D unbounded solid because frame-invariance is trivially satisfied for the only possible rigid-body motion that is, translation, in this case \cite{patnaik2020towards}. The nature of the dispersion relations obtained for the exponential and power-law kernels presented above lead to the following observations and conclusions:
\begin{itemize}
    \item The use of an exponential kernel is suitable for modeling nonlocal structures which are dispersive but not attenuating in nature. Examples of such structures can include micro- and nano-structures such as carbon nanotubes, thin films, monolayer graphene sheets, and nanobeams made of metals \cite{eringen1972linear, pradhan2009small, wang2011mechanisms, askes2011gradient, lim2015higher}. Additionally, macrostructures made of functionally graded materials and sandwiched cores, also exhibit similar dispersive characteristics \cite{romanoff2014experimental, romanoff2020review}. 
    
    \item The power-law kernel is suitable for materials which exhibit anomalous attenuation and dispersion. Note that the real and imaginary parts of the dispersion relation in Eq.~(\ref{eq: dispersion_pl}) correspond to the propagating and attenuating component of the wave. Examples of such structures can include macrostructures such as highly scattering media (e.g. fractal, porous, and layered media) and even animal tissues  \cite{fellah2004ultrasonic, fellah2004verification, magin2010fractional, buonocore2019occurrence}. Periodic media \cite{hollkamp2019analysis, hollkamp2020application, patnaik2020generalized} and structures with intentionally designed nonlocal geometric features (such as, for example, acoustic black hole metamaterials \cite{nair2019nonlocal} and metasurfaces with intentionally designed nonlocality \cite{zhu2020nonlocal}) also exhibit similar characteristics. Remarkably, the dispersion as well as the attenuation in the wave speeds in several classes of these structures, exhibit a power-law dependence on the wave-number as also captured by the use of the power-law kernel \cite{fellah2004ultrasonic, fellah2004verification, magin2010fractional, patnaik2020generalized}.
    
    \item The dispersion characteristics obtained using the exponential kernel matches exactly the dispersion relation obtained via classical differential models recovered from Eringen's approach by using exponential kernels \cite{askes2011gradient, wang2011mechanisms}. The motivation behind adopting the differential model follows from a good match of the predicted dispersion against experimental observations, particularly for different nanostructures such as, for example, graphene tubules and carbon nanotubes. As established in \cite{romano2017constitutive}, this procedure of obtaining a simplified differential model (based on the classical strain-driven approach) leads to ill-posed mathematical formulations with inconsistent predictions for different static loading conditions. On the contrary, the displacement-driven approach enables a well-posed static framework while recovering the same dynamic characteristics. 
    
    \item Since the differential approaches to nonlocal elasticity are recovered by assuming (only) the exponential kernel within Eringen's approach \cite{pradhan2009small, wang2011mechanisms, askes2011gradient, lal2019thermomechanical}, these approaches are not suitable to model structures which exhibit anomalous attenuation-dispersion characteristics. This observation suggests that the displacement-driven approach is more general, since the kernel can be modified to capture the behaviour of a wider class of size-dependent structures.
\end{itemize}
The above discussion highlights the relevance of the generalized procedure, involving analysis of both singular (e.g. power-law) and non-singular (e.g. exponential) kernels, adopted in this study. The discussion also suggests that the application of any nonlocal approach to modeling of practical structures must involve the \textit{a priori} selection of a suitable kernel, guided by experimentally observed dispersion characteristics.

\section{Application to nonlocal beams and plates}
\label{sec: application}
In this section, we will use the displacement-driven approach to analyze the static response of slender nonlocal structures, including a Timoshenko beam and a Mindlin plate. Similar to previous analyses, we will consider separately the effect of either the exponential or the power-law kernel, on the static response of the nonlocal structures subject to different loading conditions. For each structure, we will first present the nonlocal constitutive formulation, followed by the numerical results generated using the fractional-order finite element method (f-FEM). The f-FEM is a highly general formulation that is well-suited for both singular and non-singular kernels \cite{patnaik2020ritz}. The algorithm of the f-FEM is discussed in detail in \cite{patnaik2020plates} and, for the sake of brevity, details are not repeated here. For completeness, we note that the nonlocal formulation involving the exponential kernel can also be simulated by adapting the classical nonlocal finite element procedure \cite{polizzotto2001nonlocal, ansari2018bending, sidhardh2018effect}, from a strain-driven approach to a displacement-driven approach. Since we have considered shear-deformable structures in the following analyses, we adopted a reduced-order Gauss quadrature integration for the evaluation of the shear stiffness matrices in order to avoid shear-locking effects \cite{reddy2003mechanics}. More specifically, for the nonlocal Timoshenko beam we used $2$ and $1$ Gauss quadrature elements to numerically integrate the bending and shear stiffness matrices, respectively. Analogously, for the nonlocal Mindlin plate, we used $2\times 2$ and $1\times 1$ Gauss quadrature points for the bending and shear stiffness matrices, respectively. 

Analogous to classical finite element procedures, the numerical solution in the f-FEM is obtained by discretization of the Hamiltonian of the nonlocal structure by using an isoparametric formulation. Consequently, we only provide the weak-form of the governing equations for both the nonlocal Timoshenko beam and the nonlocal Mindlin plate. The corresponding strong-form of the governing equations can be easily obtained by following the derivation of the 3D governing equations in \S\ref{sec: GDE}. 
In fact, the strong-form nonlocal governing equations do not generally admit closed-from analytical expressions for the most general loading conditions. Hence, for the sake of brevity, we do not provide again the strong-form of the governing equations.

\subsection{Application to nonlocal beams}
\label{sec: beams}
\begin{figure}[h]
    \centering
    \includegraphics[width=0.5\textwidth]{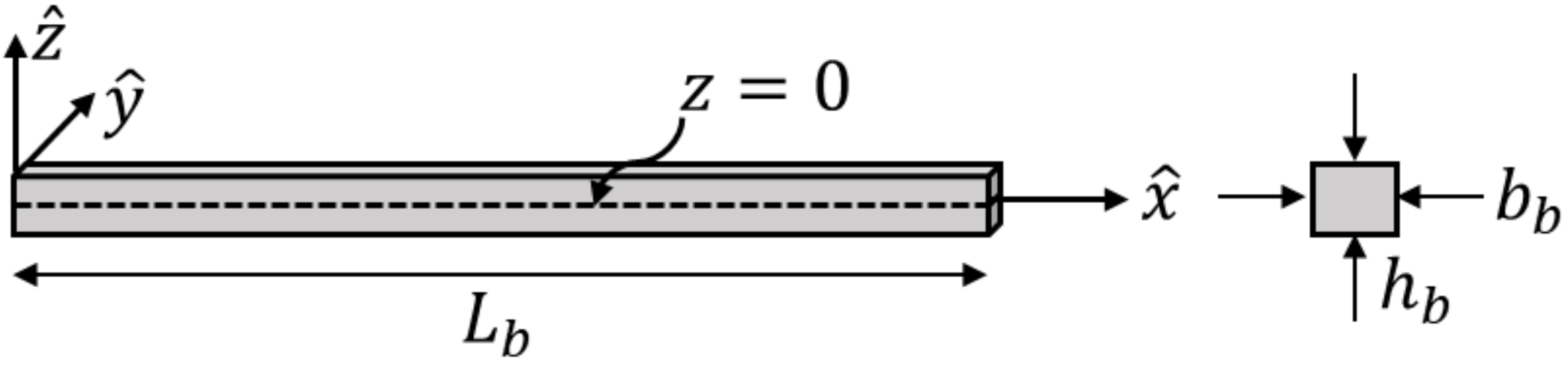}
    \caption{Schematic of the rectangular beam illustrating the different geometric parameters.}
    \label{fig: beam}
\end{figure}

In this section, we apply the displacement-driven approach to model and analyze the static response of nonlocal shear-deformable beams. We consider a nonlocal beam of uniform cross-section with length $L_b$, width $b_b$ and thickness $h_b$. A schematic of the beam, along with the chosen Cartesian coordinate axes, is provided in Fig.~(\ref{fig: beam}). As evident from Fig.~(\ref{fig: beam}), the coordinate axes are chosen such that the mid-plane of the beam is identified as $z=0$. The displacement field components at a spatial point $\bm{x}(x,y,z)$ of the beam can be expressed by following the shear-deformable Timoshenko formulation as \cite{reddy2003mechanics}:
\begin{subequations}
\label{eq: timo_disp}
\begin{equation}
    u(\bm{x}) = u_0(x) - z\theta_x(x)
\end{equation}
\begin{equation}
   w(\bm{x}) = w_0(x)
\end{equation}
\end{subequations}
where $u_0$ and $w_0$ denote the axial and transverse displacement of the mid-plane of the beam. $\theta_x$ denotes the rotation of the normal to the mid-plane of the beam, about the $\hat{x}$ axis. 

Employing the nonlocal strain-displacement relation given in Eq.~(\ref{eq: infinitesimal_fractional_strain}), we obtain the following expressions for the strains developed in the nonlocal beam:
\begin{subequations}
\label{eq: strains_timo}
\begin{equation}
    {\varepsilon}_{xx}(\bm{x}) = \overline{D}_{x}u_0(x) - z\overline{D}_{x}\theta_x(x)
\end{equation}
\begin{equation}
    {\gamma}_{xz} = \overline{D}_{x}w_0(x) - \theta_x(x)
\end{equation}
\end{subequations}
The stress developed in the beam can be be derived from the stress-strain constitutive relation in Eq.~(\ref{eq: stress}). 

Using the nonlocal strains and stresses generated in the beam, the total strain energy can be expressed as:
\begin{equation}
    \label{eq: strain_energy_timoshenko}
    \Pi = \int_0^{L_b} \int_{\frac{-b_b}{2}}^{\frac{b_b}{2}} \int_{\frac{-h_b}{2}}^{\frac{h_b}{2}} \left[ \sigma_{xx} \varepsilon_{xx} + \kappa_s \sigma_{xz} \gamma_{xz} \right] \mathrm{d}z \mathrm{d}y \mathrm{d}x
\end{equation}
where $k_s$ is the shear correction factor chosen as $\kappa_s = 5/6$ throughout this study \cite{reddy2003mechanics}. Finally the work done by external loads on the beam can be expressed as:
\begin{equation}
\label{eq: timoshenko_work}
    V = \int_0^{L_b} \int_{\frac{-b_b}{2}}^{\frac{b_b}{2}} \int_{\frac{-h_b}{2}}^{\frac{h_b}{2}} \left[ F_x u_0 + F_z w_0 + M_{\theta_x} \theta_x \right] \mathrm{d}z \mathrm{d}y \mathrm{d}x
\end{equation}
where $\{F_x,F_z\}$ are the loads applied externally in the $\hat{x}$ and $\hat{z}$ directions, respectively. Further, $M_{\theta_x}$ is the moment applied about the $\hat{y}$ axis. The nonlocal equilibrium differential equations describing the static response of the nonlocal beam can now be derived by applying variational principles ($\delta (\Pi-V)=0$).  

\subsubsection{Results and discussion}
In this section, we analyze the effect of different kernels and different loading conditions on the static bending response of the nonlocal beam. In all simulations, we considered an isotropic beam having Young's modulus $E=30$GPa and Poisson's ratio $\nu = 0.3$. The geometric dimensions of the beam were taken as $L_b=1$m and $b_b=h_b=L_b/10$. As mentioned previously, we analyzed the effect of two different kernels: the exponential kernel and the power-law kernel, on the response of the nonlocal beam. In each case, we have assumed an isotropic and symmetric horizon of nonlocality such that $l_{-_x} = l_{+_x} = l_f$ for points sufficiently far from the boundaries. Recall from \S\ref{sec: reformulation} and \S\ref{sec: FI} that this symmetry is broken ($l_{-_x}\neq l_{+_x}$) for points close to the boundaries (see Fig.~(\ref{fig: FCM})). Note that while the choice of the different material properties and geometric parameters was somewhat arbitrary, their numerical values do not affect the applicability of the model and the generality of the results presented here below.

\begin{figure}[h]
    \centering
    \includegraphics[width=\textwidth]{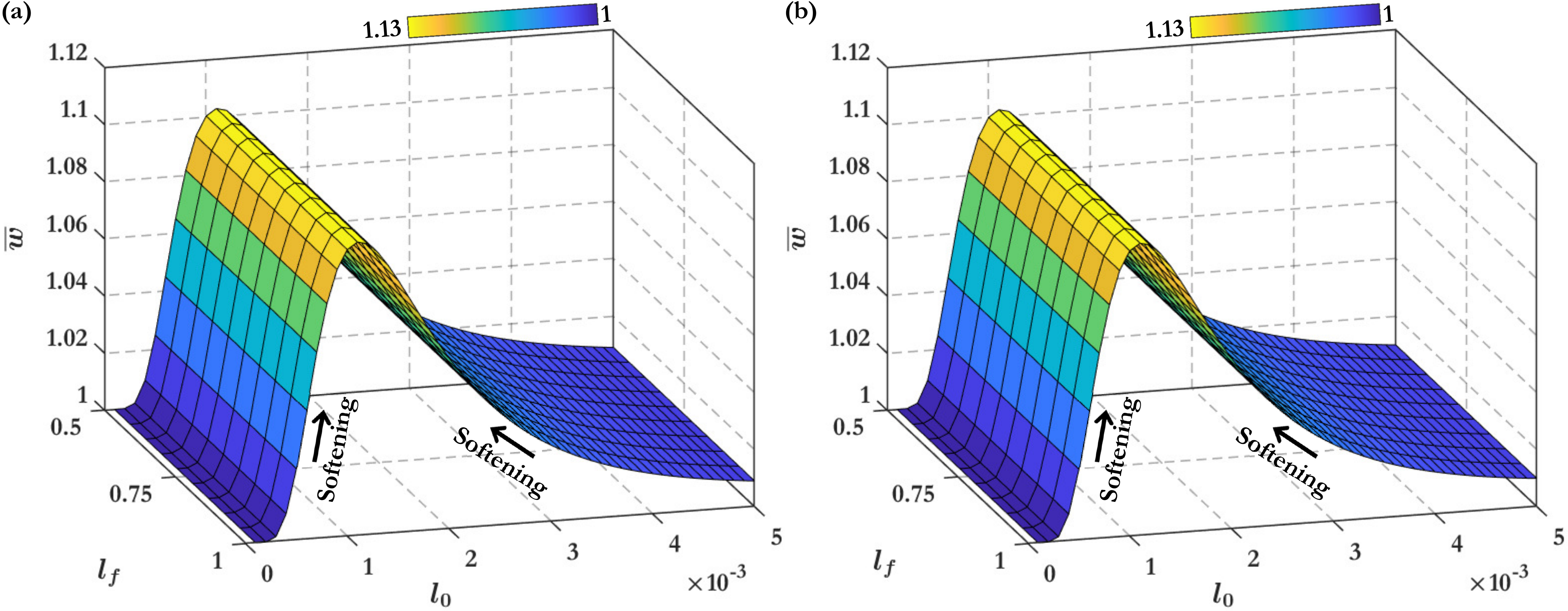}
    \caption{Bending response measured in terms of the maximum non-dimensionalized displacement of the nonlocal beam, modeled using the exponential kernel. The beam is (a) a cantilever beam subject to an transverse load at its end-point, and (b) simply-supported at all its edges and subject to a UDTL. The response is parameterized for different values of the kernel parameter $l_0$ and the nonlocal horizon $l_f$.}
    \label{fig: Exponential_beam}
\end{figure}

\begin{figure}[h]
    \centering
    \includegraphics[width=\textwidth]{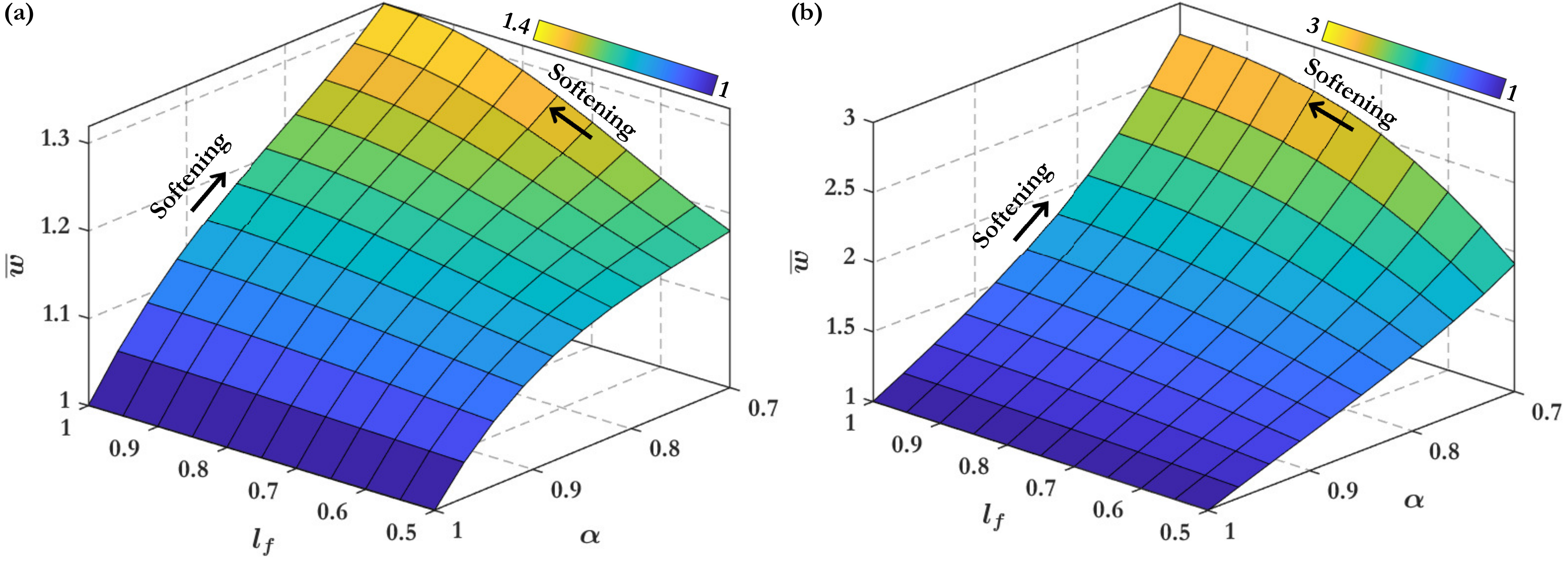}
    \caption{Bending response measured in terms of the maximum non-dimensionalized displacement of the nonlocal beam, modeled using the power-law kernel. The beam is (a) a cantilever beam subject to an transverse load at its end-point, and (b) simply-supported at all its edges and subject to a UDTL. The response is parameterized for different values of the kernel parameter $\alpha$ and the nonlocal horizon $l_f$.}
    \label{fig: frac_beam}
\end{figure}

Using the above described numerical scheme, we analyzed two different loading conditions: (a) a cantilever beam subject to a transverse force at its end-point, and (b) a simply-supported beam subject to a uniformly distributed transverse load (UDTL). For each boundary condition, we obtained the response of the beam for the following different kernel parameters:
\begin{itemize}
    \item \textit{Exponential kernel}: the kernel parameter $l_0$ was varied in $(0,0.005]\text{m}~(=(0,L_b/200])$ and the nonlocal horizon parameter $l_f$ was varied in $[0.5,1]\text{m}~ (=[L_b/2,L_b])$. The lowest value for $l_0$ was chosen as $l_0=10^{-6}$ instead of zero.
    
    \item \textit{Power-law kernel}: the nonlocal horizon $l_f$ was varied in $[0.5,1]\text{m}~ (=L_b/2,L_b])$ and the kernel parameter $\alpha$ was varied in $[0.7,1]$.
\end{itemize}
The numerical results, expressed in terms of the maximum transverse displacement of the beam, are presented in Fig.~(\ref{fig: Exponential_beam}) and Fig.~(\ref{fig: frac_beam}) for the exponential and the power-law kernel, respectively. Note that the maximum transverse displacement of the nonlocal beam is non-dimensionalized against the analogous value obtained for a local beam. The non-dimensionalized maximum transverse displacement is indicated by $\overline{w}$. Results show that the proposed displacement-driven approach predicts a consistent softening response, that is $\overline{w} > 1$, for different combinations of the kernel parameters and loading conditions. Given this general observation, we proceed to critically analyze the effects of the two different kernels and of their parameters.

The results presented for the exponential kernel in Fig.~(\ref{fig: Exponential_beam}), lead to the following observations:
\begin{itemize}
   \item \textit{\textbf{Effect of $\bm{l_0}$:}} In the limiting case $l_0 \to 0$, the classical (local) solution $\overline{w} \to 1$ is recovered \cite{romano2017nonlocal, pisano2021integral}. This observation is consistent with the fact that $l_0=0$ reduces the exponential kernel to the Dirac-delta operator, which results in a reduction of the nonlocal model to its local elastic counterpart. Note also that an increase in the value of $l_0$ first leads to an increase in the degree of softening followed by an asymptotic decline, as a function of $l_0$. Nonetheless, irrespective of the specific value of $l_0$, the nonlocal beam exhibits a softening effect, that is, $\overline{w} > 1$. This behavior is observed independently of the loading conditions and is a direct result of the nature of the exponential kernel. 
   
   Note that the strength of the nonlocal interactions between two material points, separated by a certain distance (say $d_0$) within the nonlocal beam, is directly related to the value of the nonlocal kernel. In this case, the strength of the nonlocal interactions can be given as $\propto [\exp(-d/l_0)]/l_0$ (the $l_0$ in the denominator appears from the multipliers $c^*_{-_x}$ and $c^*_{-_x}$). Consider the following function for a fixed value of $d_0$:
   \begin{equation}
       g(l_0) = \frac{1}{l_0}\exp\left(-\frac{d_0}{l_0}\right)
   \end{equation}
   A quick analysis of the above function reveals that $g(l_0)$ increases initially with an increase in $l_0$, to reach a global maximum at $l_0 = d_0$. Following the maximum, it asymptotically approaches zero with increase in $l_0$, that is, $\lim_{l_0 \to \infty} g(l_0) = 0$. It is exactly this behavior that is reflected in the static response curves of the nonlocal beam in Fig.~(\ref{fig: Exponential_beam}). Notably, this behavior is also observed in the response of atomic lattice models \cite{eringen1983differential, lim2015higher, pisano2021integral}. This observation is also consistent with our discussion in \S\ref{sec: dispersion}, where we indicated that the dispersion relation obtained from the exponential kernel is well suited to model the dynamic response of nano- and micro-structures, typically presented in the form of atomic lattice chains. Further, the result presented above is also consistent with predictions made for nonlocal beams in recently developed well-posed integral approaches \cite{pisano2021integral}.

   \item \textit{\textbf{Effect of $\bm{l_f}$:}} It appears that the length of the horizon of nonlocality $l_f$ does not significantly affect the degree of softening when $\mathcal{O}(l_f/L_b) = 1$. Recall that the parameter $l_f$ physically represents the size of the horizon of nonlocality, that is it determines the distance beyond which two particles no longer interact via long-range forces. In this context, the apparent insensitivity of $\overline{w}$ to changes in $l_f$ (in the chosen range) is a direct result of the fact that $\mathcal{K}(x,x^\prime)=e^{-|x-x^\prime|/l_0} \approx 0$ for $|x-x^\prime| > 5l_0$. Thus, the effect of the nonlocal kernel is always accounted completely (in the chosen range for $l_f$) and does not change appreciably with a further increase in $l_f$. However, note that $l_f$ plays a critical role (along with $c^*_{-_x}$ and $c^*_{+_x}$) in achieving frame-invariance and consistency to different boundary conditions.
\end{itemize}

We now proceed to analyze the results obtained for the nonlocal beam modeled via the power-law kernel. A detailed analysis of the results presented in Fig.~(\ref{fig: frac_beam}) leads to the following conclusions:
\begin{itemize}
    \item \textit{\textbf{Effect of $\bm\alpha$:}} For $\alpha=1$ we recover the classical (local) solution $\overline{w} = 1$, consistent with the nature of the power-law kernel \cite{patnaik2020generalized,sumelka2014thermoelasticity}. Next, observe that a decrease in the value of $\alpha$ leads to a consistent increase in the degree of softening, irrespective of the nature of the loading conditions. A decrease in the value of $\alpha$ leads to an increase in the value of the power-law kernel $|x-x^\prime|^{-\alpha}$ for a fixed value of $|x-x^\prime|$. This suggests that a decrease in the value of $\alpha$ leads to an increase in the degree of nonlocality and a higher degree of softening. We merely note that the above findings are also consistent with results presented in existing fractional-order approaches to nonlocal elasticity (which employ the power-law kernel) \cite{patnaik2020ritz,patnaik2020plates}. We highlight that there exists a critical value of $\alpha$ ($=0.4$ for the nonlocal beam), beyond which the beam undergoes excessive softening and the model breaks down. More detailed discussions on this aspect can be found in \cite{patnaik2020ritz, sidhardh2020geometrically, sumelka2014thermoelasticity}.
    
    \item \textit{\textbf{Effect of $\bm{l_f}$:}} An increase in the value of $l_f$ leads to significant softening in the case of a power-law kernel. The enhanced softening is a direct result of the fact that, an increase in the value of $l_f$ leads to an increase in the size of the horizon of nonlocality. Consequently, a larger number of points within the solid contribute to the nonlocal interactions leading to an increase in the degree of nonlocality. Again, this observation is supported by results presented in fractional-order approaches to nonlocal elasticity \cite{patnaik2020ritz,patnaik2020plates}. Note that, unlike constraints on the lower limit for $\alpha$, the value of $l_f$ can approach zero. In fact, in the limiting case $l_f \to 0$, we recover the classical local solution \cite{sidhardh2020thermodynamics,patnaik2020towards}.
\end{itemize}

\subsection{Application to nonlocal plates}
\label{sec: plates}
\begin{figure}[h]
    \centering
    \includegraphics[width=0.4\textwidth]{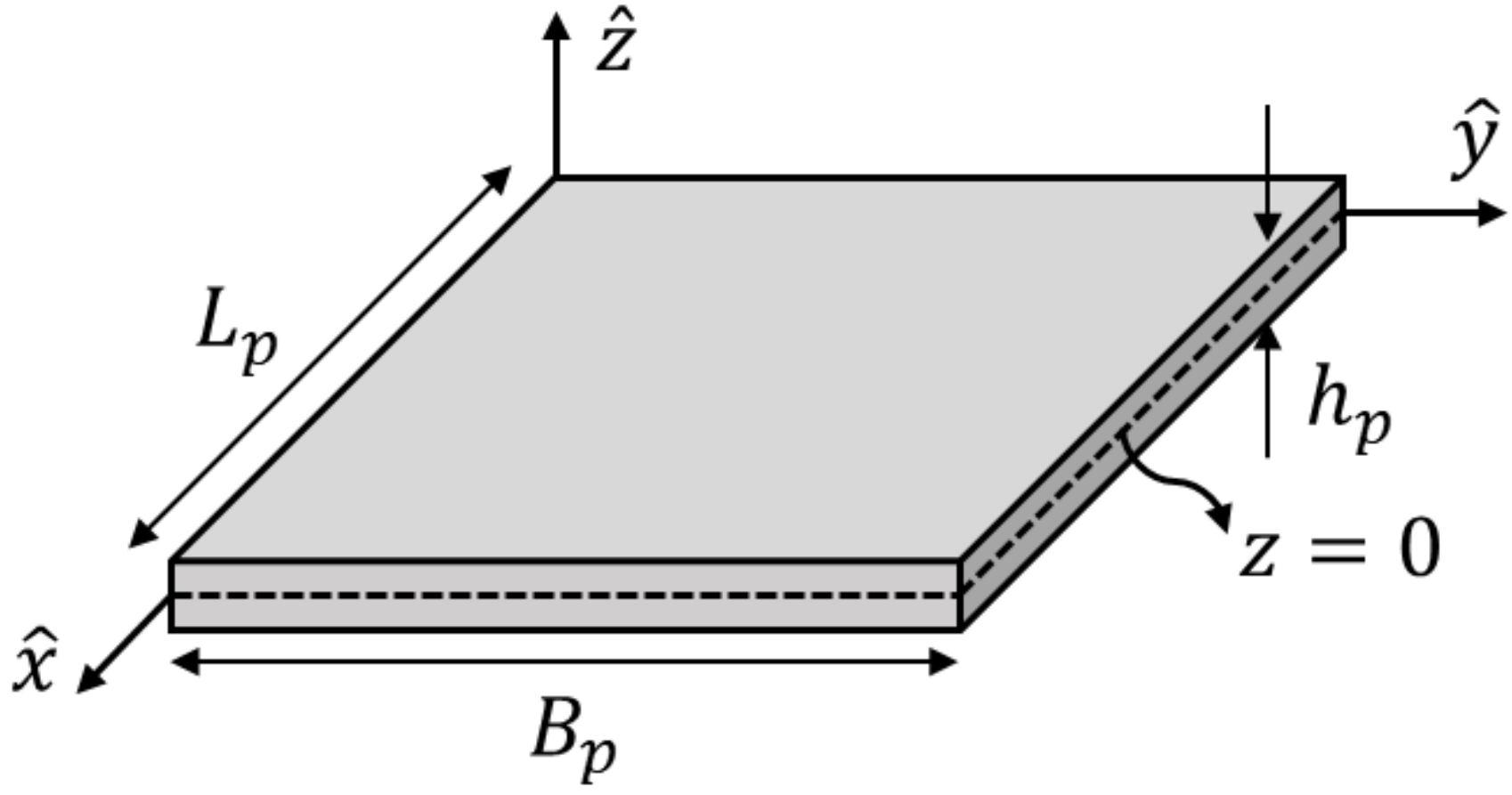}
    \caption{Schematic of the rectangular plate illustrating the different geometric parameters.}
    \label{fig: plate}
\end{figure}

We conclude this study by analyzing the response of nonlocal shear-deformable plate, in a manner similar to the analysis of nonlocal beams in \S\ref{sec: beams}. For this purpose, we consider a rectangular plate of uniform thickness, illustrated in Fig.~(\ref{fig: plate}). The length, width, and thickness of the plate are denoted by $L_p$, $B_p$ and $h_p$, respectively. The Cartesian reference frame, as indicated in Fig.~(\ref{fig: plate}), is chosen such that $z=0$ denotes the mid-plane of the plate. In the chosen reference frame, the displacement field components at a spatial point $\bm{x}(x,y,z)$ can be expressed by following the shear-deformable Mindlin formulation as \cite{reddy2003mechanics}:
\begin{subequations}
\label{eq: Mindlin_Kinematics}
    \begin{equation}
    u(\bm{x}) = u_0(x,y) - z\theta_x(x,y)
    \end{equation}
    \begin{equation}
    v(\bm{x}) = v_0(x,y) - z\theta_y(x,y)
    \end{equation}
    \begin{equation}
    \label{eq: tranverse_displacment}
    w(\bm{x}) = w_0(x,y)
    \end{equation}
\end{subequations}
where $u_0$, $v_0$, and $w_0$ denote the displacements at the mid-plane of the plate along the $\hat{x}$, $\hat{y}$, and $\hat{z}$ directions. $\theta_x$ and $\theta_y$ denote the rotations of the transverse normal at the mid-plane of the plate, about the $\hat{y}$ and $\hat{x}$ axes. 

The nonzero infinitesimal nonlocal strains, corresponding to the deformation field in Eq.~(\ref{eq: Mindlin_Kinematics}), are obtained by using Eq.~(\ref{eq: infinitesimal_fractional_strain}) as:
\begin{subequations}
\label{eq: Mindlin_strains}
\begin{equation}
    \varepsilon_{xx}(\bm{x}) = \overline{D}_x u_0(x,y) - z \overline{D}_x \theta_x(x,y)
\end{equation}
\begin{equation}
    \varepsilon_{yy}(\bm{x}) = \overline{D}_y v_0(x,y) - z \overline{D}_y \theta_y(x,y)
\end{equation}
\begin{equation}
    \gamma_{xy}(\bm{x}) = \overline{D}_y u_0(x,y) + \overline{D}_x v_0(x,y) - z \left[ \overline{D}_y \theta_x(x,y) + \overline{D}_y \theta_y(x,y) \right]
\end{equation}
\begin{equation}
    \gamma_{xz}(\bm{x}) = \overline{D}_x w_0(x,y) - \theta_x(x,y)
\end{equation}
\begin{equation}
    \gamma_{yz}(\bm{x}) = \overline{D}_y w_0(x,y) - \theta_y(x,y)
\end{equation}
\end{subequations}
The nonlocal stresses developed in the nonlocal plate can now be evaluated using the stress-strain constitutive relation in Eq.~(\ref{eq: stress}). 

Using the nonlocal strains and stresses generated in the plate, the total strain energy can be expressed as:
\begin{equation}
    \label{eq: strain_energy_mindlin}
    \Pi = \int_0^{L_p} \int_0^{B_p} \int_{\frac{-h_p}{2}}^{\frac{h_p}{2}} \left[ \sigma_{xx} \varepsilon_{xx} + \sigma_{yy} \varepsilon_{yy} + \sigma_{xy} \gamma_{xy} + \kappa_s \sigma_{xz} \gamma_{xz} + \kappa_s \sigma_{yz} \gamma_{yz} \right] \mathrm{d}z \mathrm{d}y \mathrm{d}x
\end{equation}
Finally the work done by external loads on the plate can be expressed as:
\begin{equation}
\label{eq: mindlin_work}
    V = \int_0^{L_p} \int_0^{B_p} \int_{\frac{-h_p}{2}}^{\frac{h_p}{2}}  \left[ F_x u_0 + F_y v_0 + F_z w_0 + M_{\theta_x} \theta_x + M_{\theta_y} \theta_y \right] \mathrm{d}z \mathrm{d}y \mathrm{d}x
\end{equation}
where $\{F_x,F_y,F_z\}$ are the loads applied externally in the $\hat{x}$, $\hat{y}$, and $\hat{z}$ directions, respectively. Further, $\{M_{\theta_x},M_{\theta_y}\}$ are the moments applied about the $\hat{y}$ and $\hat{x}$ axes, respectively. The nonlocal equilibrium differential equations describing the static response of the nonlocal plate can now be derived by applying variational principles. The response of the nonlocal plate is numerically simulated using the f-FEM and the results are presented in the following section.

\subsubsection{Results and discussion}
\begin{figure}[h]
    \centering
    \includegraphics[width=\textwidth]{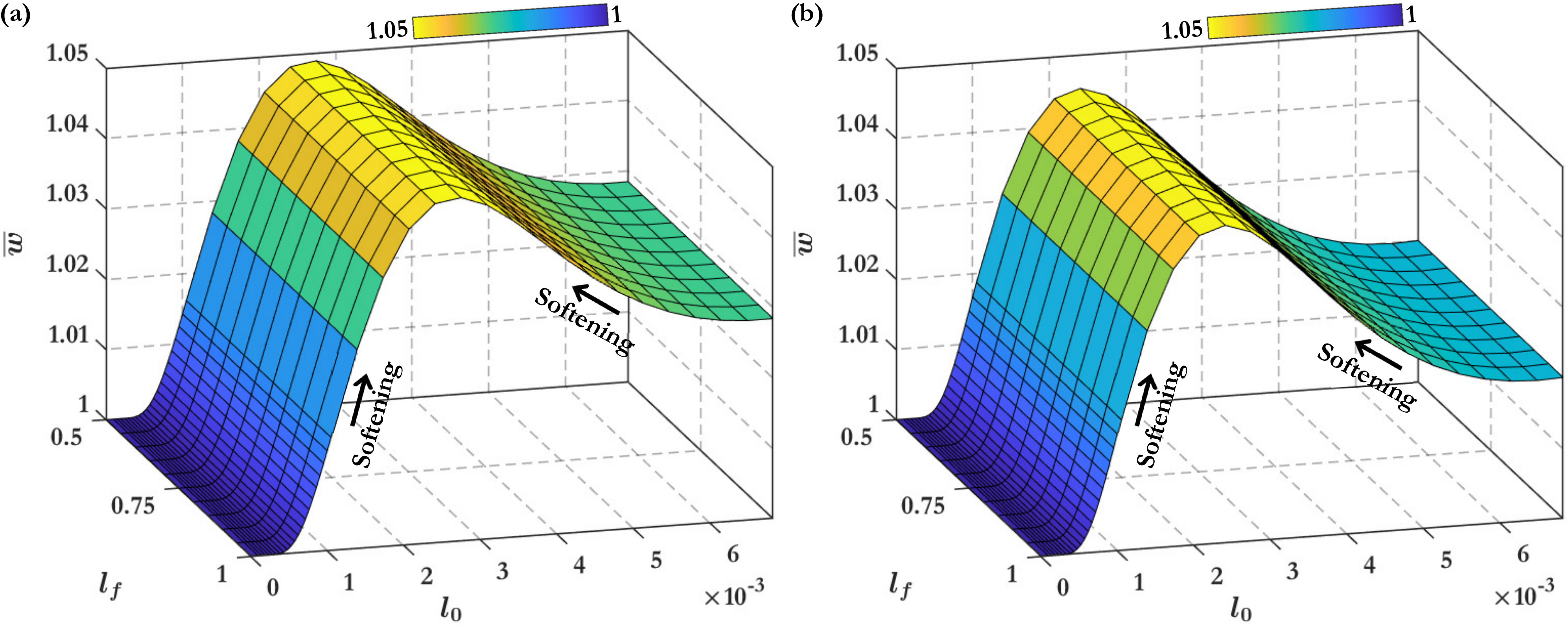}
    \caption{Bending response measured in terms of the non-dimensionalized displacement at the center point of the nonlocal plate subject to a UDTL and modeled using the exponential kernel. The plate is (a) clamped at all its edges and (b) simply-supported at all its edges. The response is parameterized for different values of the kernel parameter $l_0$ and the nonlocal horizon $l_f$.}
    \label{fig: Exponential}
\end{figure}

\begin{figure}[h!]
    \centering
    \includegraphics[width=\textwidth]{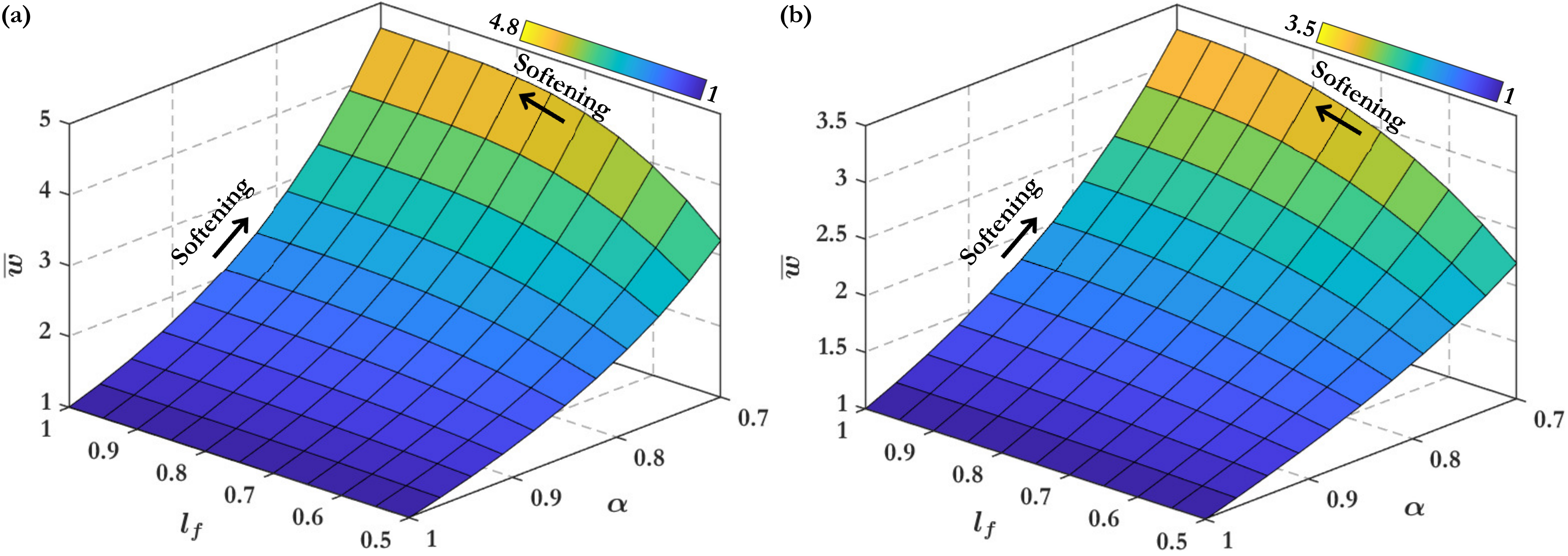}
    \caption{Bending response measured in terms of the non-dimensionalized displacement at the center point of the nonlocal plate subject to a UDTL and modeled using the power-law kernel. The plate is (a) clamped at all its edges and (b) simply-supported at all its edges. The response is parameterized for different values of the kernel parameter $\alpha$ and the nonlocal horizon $l_f$.}
    \label{fig: frac}
\end{figure}

In this section, we analyze the effect of different kernels and different boundary conditions on the static bending response of the nonlocal plate. In all the simulations, we considered an isotropic plate with material properties identical to the nonlocal beam in \S\ref{sec: beams}. The in-plane dimensions of the plate were chosen as $L_p=B_p=1$m, and thickness was chosen as $h_p=L_p/10$. We considered the effect of the exponential kernel and the power-law kernel, similar to the analysis conducted for the nonlocal beam. In each case, we have assumed an isotropic and symmetric horizon of nonlocality such that $l_{-_j} = l_{+_j} = l_f$ for points sufficiently far from the boundaries. We analyzed the static response of the plate subject to a uniformly distributed transverse load (UDTL) for two different kinds of boundary conditions: (a) the plate clamped at all the edges, and (b) the plate simply-supported at all its edges, for different combinations of the kernel parameters. 
The constraints on the generalized displacement coordinates corresponding to these boundary conditions are as follows \cite{reddy2003mechanics}:
\begin{subequations}
\label{eq: bcs}
\begin{equation}
\label{eq: CCCC}
\text{Clamped plate}: \left\{
\begin{matrix}
& \hat{x} = \{0,L_p\} : u_0 = v_0 = w_0 = \theta_x = \theta_y = 0\\
& \hat{y} = \{0,L_p\} : u_0 = v_0 = w_0 = \theta_x = \theta_y = 0
\end{matrix}
\right.
\end{equation}
\begin{equation}
\label{eq: SSSS}
\text{Simply-supported plate}: \left\{
\begin{matrix}
& \hat{x} = \{0,L_p\} : v_0 = w_0 = \theta_y = 0\\
& \hat{y} = \{0,B_p\} : u_0 = w_0 = \theta_x = 0
\end{matrix}
\right.
\end{equation}
\end{subequations}
For each boundary condition, we obtained the response of the plate for the following different kernel parameters:
\begin{itemize}
    \item \textit{Exponential kernel}: the kernel parameter $l_0$ was varied in $(0,0.005]\text{m}~(=(0,L_p/200])$ and the nonlocal horizon parameter $l_f$ was varied in $[0.5,1]\text{m}~ (=[L_p/2,L_p])$. The lowest value for $l_0$ was chosen as $l_0=10^{-6}$ instead of zero. This is because the choice of $l_0=0$ leads to a singularity within the exponential kernel, wherein the exponential kernel reduces to a Dirac-delta operator \cite{romano2017constitutive}.
    
    \item \textit{Power-law kernel}: the nonlocal horizon parameter $l_f$ was varied in $[0.5,1]\text{m}~ (=[L_p/2,L_p])$ and the kernel parameter $\alpha$, also called as the fractional-order \cite{patnaik2020generalized}, was varied in $[0.7,1]$. Similar to the nonlocal beam, there exists a constraint on the lower limit of $\alpha$ (=$0.4$ for the nonlocal plate), while the value of $l_f$ can approach to zero \cite{patnaik2020plates}.
\end{itemize}
The numerical results, in terms of the maximum transverse displacement obtained at the mid-point of the plate, are presented in Fig.~(\ref{fig: Exponential}) and Fig.~(\ref{fig: frac}) for the exponential and the power-law kernel, respectively. Analogous to the nonlocal beam, the transverse displacement of the mid-point of the nonlocal plate is non-dimensionalized against the analogous value obtained for a local Mindlin plate. As evident from the results, the same considerations noted for the nonlocal Timoshenko beam directly extend to the nonlocal Mindlin plate.

\subsection{Few additional remarks}
The displacement-driven approach predicts a consistent softening response for both nonlocal beams and plates, irrespective of the nature of the kernel and the loading conditions. This is a significant outcome when compared to the often incoherent predictions of classical (integer-order) strain-driven integral approaches \cite{challamel2008small,fernandez2016bending,romano2017constitutive}. The results obtained from the displacement-driven approach are free from stiffening effects or even the absence of nonlocal effects noted in some strain-driven methods, under selected loading conditions. This latter statement does not suggest that stiffening effects are not expected in nonlocal structures; it only implies that stiffening effects should not be predicted by strain-driven approaches \cite{romano2017constitutive, pisano2021integral}.

Indeed, several nonlocal structures do exhibit stiffening effects \cite{romano2017constitutive, patnaik2020towards, pisano2021integral}. In this context, we highlight here that the displacement-driven approach can be extended in a straight-forward and physically consistent manner to capture both stiffening and softening effects. This has been established exclusively for power-law kernels in \cite{patnaik2020towards}. Further, the displacement-driven approach is also well-suited for the geometrically nonlinear analyses of nonlocal structures. The extension to geometrically nonlinear analyses has been performed in \cite{sidhardh2020geometrically, patnaik2020geometrically, sidhardh2020analysis} to analyse the nonlinear bending and postbuckling response of nonlocal beams and plates described by power-law kernels. The above discussion suggests that the displacement-driven approach presents a promising alternative to the analysis of nonlocal structures.

\section{Conclusions}
We presented an approach to model the nonlocal behavior in an elastic solid that leverages strain-displacement relations formulated in integral form. This method was dubbed displacement-driven approach to clearly differentiate it from the existing stress- and strain-driven techniques. The mathematical, physical, and thermodynamic consistency of the method were carefully demonstrated and shown to address some critical deficiencies of existing integral-type models of nonlocal elasticity. 
This result may also be seen as a generalization of fractional order formulations of nonlocal elasticity that restricts the nonlocal kernel to power-law type.
The well-posedness of the elastic governing equations, which is a requisite for modeling well-posed elastic behavior, is achieved by the present approach. This observation is a direct consequence of the positive-definite and convex deformation energy density achievable by using the displacement-driven approach. The formulation also presents some other important benefits including the generalization to include asymmetric kernels arising from material heterogeneities, and the proper truncation of the nonlocal horizon at physical boundaries to maintain the completeness of the kernel. These features will be fundamental as nonlocal theories start percolating into a broad spectrum of practical applications. We also derived the dispersion relations for a 1D infinite solid and showed the practical influence of some commonly used nonlocal kernels. It was found that different kernels capture different features of the dynamic behavior, hence their selection for dynamic analyses should be properly pondered.
Finally, numerical simulations were performed to investigate the static response of both nonlocal Timoshenko beams and nonlocal Mindlin plates formulated based on the displacement-driven approach. Numerical results highlighted the very consistent physical behavior enabled by the displacement-driven approach. This consistency was observable based on the monotonic softening of the elastic response following an increase in the degree of nonlocality; a behavior observed irrespective of the choice of the nonlocal kernel, and of the loading and boundary conditions. In conclusion, the displacement-driven formulation provides a powerful framework to model nonlocal effects in elastic solids while enabling the use of different kernels that are representative of different dynamic material behavior.

\section*{Acknowledgements}
The authors gratefully acknowledge the financial support of the National Science Foundation (NSF) under grants MOMS \#1761423 and DCSD \#1825837, and the Defense Advanced Research Project Agency (DARPA) under grant \#D19AP00052. S.P. acknowledges the partial support received through the R. H. Kohr Graduate Fellowship. The content and information presented in this manuscript do not necessarily reflect the position or the policy of the government. The material is approved for public release; distribution is unlimited.

\section*{Competing Interests} 
The authors declare that they have no conflict of interest.

\section*{Appendix: Derivation of the strong-form of equilibrium equations}
\label{sec: Appendix_B}
In the following, we present the details of the steps involved in the derivation of differ-integral (nonlocal) equilibrium equation for the nonlocal solid. For this purpose, we recall the Hamilton's principle in Eq.~(\ref{eq: extended_hamiltons_principle2}).
\begin{equation}
    \tag{\ref{eq: extended_hamiltons_principle2}}
    \int_{t_1}^{t_2} \delta\left[ \underbrace{ \frac{1}{2} \int_\Omega (\bm\sigma:\bm\varepsilon) \mathrm{d}\mathbb{V}}_{\text{Strain Energy}} - \underbrace{\int_\Omega (\bm{\overline{b}}\cdot\bm{u}) \mathrm{d}\mathbb{V}  - \int_{\partial \Omega} (\bm{\overline{t}}\cdot\bm{u}) \mathrm{d}\mathbb{A}}_{\text{External Work}} -  \underbrace{\frac{1}{2} \int_{\Omega} \rho (\bm{\dot{u}} \cdot {\bm{\dot{u}}})\mathrm{d}\mathbb{V}}_{\text{Kinetic Energy}} \right] \mathrm{d}t=0
\end{equation}
The variations of terms corresponding to the external work and kinetic energy are identical to their counterparts in local elasticity and hence, in the following, we will derive only the variation of the strain energy:
\begin{equation}
    \label{eq: I1}
    \mathcal{I}_1=\int_{t_1}^{t_2}\delta\left[ \frac{1}{2} \int_\Omega (\bm\sigma:\bm\varepsilon) \mathrm{d}\mathbb{V}\right]\mathrm{d}t
\end{equation}
The above expression can be simplified by using the symmetries in the strain and stress tensors as:
\begin{equation}
    \label{eq: simp1}
     \mathcal{I}_1 =\int_{t_1}^{t_2}\underbrace{\int_\Omega {\sigma}_{ij}(\bm{x})~\delta \overline{D}_{x_j}\left[u_i(\bm{x})\right] \mathrm{d}\mathbb{V}}_{\mathcal{I}_2}\mathrm{d}t
\end{equation}
Now, by using the definition of the nonlocal operator $\overline{D}_{x_j}(\cdot)$ from Eq.~(\ref{eq: operator_definition_2}) we obtain:
\begin{equation}
\label{eq: stress_step2}
\begin{split}
    \mathcal{I}_2 = \int_{\Omega_{x_i}}\int_{\Omega_{x_k}} \int_{\Omega_{x_j}} \sigma_{ij}(\bm{x}) \left[ c_{-_j}^*  \int_{x_{-j}}^{x_j} \mathcal{K}(\bm{x},\bm{x}^\prime) \left[\frac{\mathrm{d}\delta u_i(x^\prime_j)}{\mathrm{d}x_j^\prime} \right] \mathrm{d}x_j^\prime + c_{+_j}^*  \int_{x_{j}}^{x_{+j}} \mathcal{K}(\bm{x},\bm{x}^\prime) \left[\frac{\mathrm{d}\delta u_i(x^\prime_j)}{\mathrm{d}x_j^\prime} \right] \mathrm{d}x_j^\prime \right] \mathrm{d}x_j \mathrm{d}x_k \mathrm{d}x_i
\end{split}
\end{equation}
Each of the two terms within the inner integral in the above expression, can be re-written by changing the order of integration over the nonlocal spatial variable $x_j^\prime$ and local spatial variable along the same direction $x_j$. Recall also that $x_{-j}=x_j-l_-$ and $x_{+j}=x_j+l_+$. This gives:
\begin{subequations}
\begin{equation}
    \mathcal{I}_{21}= c_{-_j}^* \int_{\Omega_{x_j}} \sigma_{ij}(\bm{x}) \int_{x_{j}-l_-}^{x_j} \mathcal{K}(\bm{x},\bm{x}^\prime) \left[ \frac{\mathrm{d}\delta u_i(x^\prime_j)}{\mathrm{d}x_j^\prime} \right] \mathrm{d}x_j^\prime \mathrm{d}x_j = c_{-_j}^* \int_{\Omega_{x_j}} \left[ \frac{\mathrm{d}\delta u_i(x^\prime_j)}{\mathrm{d}x_j^\prime} \right] \underbrace{\int_{x_{j}}^{x_j+l_-} \mathcal{K}(\bm{x},\bm{x}^\prime)  \sigma_{ij}(\bm{x}) \mathrm{d}x_j}_{{I}^-_{x_j^\prime} \left[ \sigma_{ij}(x_j^\prime) \right]}\mathrm{d}x_j^\prime
\end{equation}
\begin{equation}
     \mathcal{I}_{22}= c_{+_j}^* \int_{\Omega_{x_j}} \sigma_{ij}(\bm{x}) \int_{x_{j}}^{x_j+l_+} \mathcal{K}(\bm{x},\bm{x}^\prime) \left[ \frac{\mathrm{d}\delta u_i(x^\prime_j)}{\mathrm{d}x_j^\prime} \right] \mathrm{d}x_j^\prime\mathrm{d}x_j = c_{+_j}^* \int_{\Omega_{x_j}} \left[ \frac{\mathrm{d}\delta u_i(x^\prime_j)}{\mathrm{d}x_j^\prime} \right] \underbrace{\int_{x_{j}-l_+}^{x_j} \mathcal{K}(\bm{x},\bm{x}^\prime)  \sigma_{ij}(\textbf{x}) \mathrm{d}x_j}_{{I}^+_{x_j^\prime} \left[ \sigma_{ij}(x_j^\prime) \right]}\mathrm{d}x_j^\prime
\end{equation}
\end{subequations}
The indicated expressions in the above equation, corresponding to the convolution of the nonlocal stress tensor, are denoted by ${I}^-_{x_j} \left[ \sigma_{ij}(x_j^\prime) \right]$ and ${I}^+_{x_j} \left[ \sigma_{ij}(x_j^\prime) \right]$. We again use integration-by-parts to separate the variation of displacement field in the above expressions which results in the following expressions:
\begin{subequations}
\begin{equation}
\begin{split}
   \mathcal{I}_{21} = c_{-_j}^* \left[ {I}^-_{x_j} \left[ \sigma_{ij}(x_j^\prime) \right] \right] \delta u_i(x^\prime_j) \bigg\vert_{\partial \Omega_{x_j}}
    - \int_{\Omega_{x_j}} \delta u_i(x_j^\prime) \left[ \frac{\mathrm{d}}{\mathrm{d}x_j^\prime} \left( c_{-_j}^* {I}^-_{x_j} \left[ \sigma_{ij}(x_j^\prime) \right] \right) \right] \mathrm{d}x_j^\prime
\end{split}
\end{equation}
\begin{equation}
   \mathcal{I}_{22} = c_{+_j}^* \left[ {I}^+_{x_j} \left[ \sigma_{ij}(x_j^\prime) \right] \right] \delta u_i(x^\prime_j) \bigg\vert_{\partial \Omega_{x_j}}
    - \int_{\Omega_{x_j}} \delta u_i(x_j^\prime) \left[ \frac{\mathrm{d}}{\mathrm{d}x_j^\prime} \left( c_{+_j}^* {I}^+_{x_j} \left[ \sigma_{ij}(x_j^\prime) \right] \right) \right] \mathrm{d}x_j^\prime
\end{equation}
\end{subequations}
Combining the above results yields the following simplified expression:
\begin{equation}
\label{eq: stress_step3}
    \mathcal{I}_2 = -\int_{\Omega_{x_i}}\int_{\Omega_{x_k}}\int_{\Omega_{x_j}}  \left[ \frac{\mathrm{d}\tilde{I}_{x_j}\sigma_{ij}(x_j)}{\mathrm{d}x_j} \right] \delta u_i(x_j)\mathrm{d}x_j \mathrm{d}x_k \mathrm{d}x_i+\int_{\Omega_{x_i}}\int_{\Omega_{x_k}}\left\{\tilde{I}_{x_j}\sigma_{ij}(x_j)\right\}\bigg\vert_{\partial \Omega_{x_j}} \mathrm{d}x_k \mathrm{d}x_i
\end{equation}
where, $\tilde{I}_{x_j}\sigma_{ij}(x_j)$ follows the definition of the integral operator defined in Eq.~(\ref{eq: integral_def}). In obtaining the above result, we interchanged the dummy variables $x_j^\prime$ and $x_j$ are interchanges below for the sake of brevity. Finally, by substituting the above result in Eq.~(\ref{eq: simp1}), we obtain:
\begin{equation}
    \label{eq: final_stress_simplified}
    \mathcal{I}_1= \int_{t_1}^{t_2}\left[-\int_{\Omega} \left[ \left(  \bm{\widetilde{\nabla}}^{\alpha} \cdot \bm{\sigma} \right) \cdot \delta \bm{u} \right] \mathrm{d}\mathbb{V} + \\
    \int_{\partial \Omega} \left[ \left(  \bm{I}^{1-\alpha_1}_{\bm{\hat{n}}} \cdot \bm{\sigma} \right) \cdot \delta \bm{u} \right] \mathrm{d}\mathbb{A}\right]\mathrm{d}t 
\end{equation}
Applying the fundamental principle of variational calculus, we recover the governing equations in Eq.~\eqref{eq: GDE} and boundary conditions in Eq.~\eqref{eq: BCs}.

\bibliographystyle{naturemag}
\bibliography{Report}
\end{document}